\documentclass[11pt]{article}

\usepackage{amssymb}
\usepackage{amsfonts}
\usepackage{amsbsy}

\newcommand{\ra}{{\rightarrow}}
\newcommand{\eproof}{\hfill\rule{2.2mm}{3.0mm}}
\newcommand{\Proof}{\noindent {\bf Proof.~~}}
\newcommand{\ep}{\varepsilon}
\newcommand{\wmod}[1]{\mbox{~(mod~$#1$)}}
\newcommand{\shsp}{\hspace{1em}}

\newcommand{\A}{\he{A}}
\newcommand{\bH}{{\bf H}}
\renewcommand{\vec}[1]{{\bf #1}}

\newtheorem{prop}{Proposition}[section]
\newtheorem{lemma}[prop]{Lemma}
\newtheorem{defi}{Definition}[section]

\newtheorem{theorem}[prop]{Theorem}

\newtheorem{exam}{Example}[section]

\def\slfrac#1#2{\hbox{\kern.1em %
 \raise.5ex\hbox{\the\scriptfont0 #1}\kern-.11em %
 /\kern-.15em\lower.25ex\hbox{\the\scriptfont0 #2}}}

\font\phvr=phvr at 11pt
\newcommand\he[1]{\mbox{\phvr  #1}}
\newcommand{\eqn}[1]{(\ref{#1})}
\newcommand{\hsp}{\hspace*{\parindent}}

\newcommand{\eeq}{\end{equation}}
\newcommand{\beql}[1]{\begin{equation}\label{#1}}

\newcommand\us[1]{\underline{#1}}

\newcommand{\fa}{\mbox{for all}}

\newcommand{\rdet}{\rm det}

\newcommand{\la}{\lambda}

\newcommand{\ZZ}{{\mathbb Z}}
\newcommand{\RR}{{\mathbb R}}

\newcommand{\KK}{{\mathbb K}}

\newcommand{\CC}{{\mathbb C}}
\newcommand{\bo}{{\bf 0}}
\newcommand{\bt}{{\bf t}}
\newcommand{\bx}{{\bf x}}
\newcommand{\bv}{{\bf v}}

\newcommand{\bd}{{\bf d}}
\newcommand{\by}{{\bf y}}

\newcommand{\sD}{{\cal D}}
\newcommand{\sG}{{\cal G}}

\newcommand{\sS}{{\cal S}}
\newcommand{\sT}{{\cal T}}

\newcommand{\sX}{{\cal X}}

\makeatletter
\def\@sect#1#2#3#4#5#6[#7]#8{\ifnum #2>\c@secnumdepth
     \def\@svsec{}\else
     \refstepcounter{#1}\edef\@svsec{\csname the#1\endcsname.\hskip .75em }\fi
     \@tempskipa #5\relax
      \ifdim \@tempskipa>\z@
        \begingroup #6\relax
          \@hangfrom{\hskip #3\relax\@svsec}{\interlinepenalty \@M #8\par}%
        \endgroup
       \csname #1mark\endcsname{#7}\addcontentsline
         {toc}{#1}{\ifnum #2>\c@secnumdepth \else
                      \protect\numberline{\csname the#1\endcsname}\fi
                    #7}\else
        \def\@svsechd{#6\hskip #3\@svsec #8\csname #1mark\endcsname
                      {#7}\addcontentsline
                           {toc}{#1}{\ifnum #2>\c@secnumdepth \else
                             \protect\numberline{\csname the#1\endcsname}\fi
                       #7}}\fi
     \@xsect{#5}}
\def\@begintheorem#1#2{\it \trivlist \item[\hskip \labelsep{\bf #1\ #2.}]}

\def\plain{plain}\ifx\fmtname\plain\csname fi\endcsname
     
     \input docstrip
     \preamble

     Do not distribute the stripped version of this file.
     The checksum in the header refers to the documented version.

     \endpreamble
     \generateFile{here.sty}{t}{\from{here.doc}{}}
     \endinput
\fi
\ifcat a\noexpand @\let\next\relax\else\def\next{%
    \documentstyle[here,doc]{article}\MakePercentIgnore}\fi\next
\ifx\@Hxfloat\@Hundef\else\expandafter\endinput\fi
\let\@Hxfloat\@xfloat
\def\@xfloat#1[{\@ifnextchar{H}{\@HHfloat{#1}[}{\@Hxfloat{#1}[}}
\def\@HHfloat#1[H]{%
\expandafter\let\csname end#1\endcsname\end@Hfloat
\vskip\intextsep\vbox\bgroup\def\@captype{#1}\parindent\z@
\ignorespaces}
\def\end@Hfloat{\egroup\vskip \intextsep}

\makeatother

\begin{document}
\begin{center}
{\Large {\bf  Substitution Delone Sets}} \\
\vspace{1.5\baselineskip}
{\em Jeffrey C. Lagarias} \\
\vspace*{.2\baselineskip}
AT\&T Labs - Research\\
Florham Park, NJ 07932-0971 \\
\vspace{1\baselineskip}
{\em Yang Wang} \\
School of Mathematics \\
Georgia Institute of Technology \\
Atlanta, GA 30332 \\
\vspace*{1\baselineskip}
(August 31, 2002 version) \\
\vspace{1.5\baselineskip}
{\bf ABSTRACT}
\end{center}

Substitution Delone set families   are families
of Delone sets $\sX =(X_1, \dots,
X_n)$  which satisfy the inflation functional
equation
$$
X_i = \bigvee_{j=1}^m ( \he{A} (X_j) + \sD_{ij} ) , \quad
1 \le i \le m ~,
$$
in which  $\he{A}$  is an expanding matrix, i.e. all of the
eigenvalues of $\he{A}$ fall outside the unit circle.
Here the $\sD_{ij}$ are finite sets of vectors in $\RR^d$ 
and $\bigvee$ denotes union that counts multiplicity. 

This paper
characterizes families  $\sX=(X_1, ..., X_n)$
 that satisfy an inflation
functional equation, in which each $X_i$ is
a multiset (set with multiplicity) whose underlying set is discrete.
It then studies the subclass of Delone set solutions, and
gives necessary conditions on the coefficients of
the inflation functional equation for such solutions $\sX$
to exist. It relates Delone set solutions to a narrower
 subclass of solutions, called self-replicating multi-tiling sets, 
which arise as tiling
sets for self-replicating multi-tilings. \\



\noindent
\clearpage
\thispagestyle{empty}
\setcounter{page}{1}
\vspace*{1\baselineskip}
\begin{center}
{\Large {\bf  Substitution Delone Sets}} \\
\vspace{1.5\baselineskip}
{\em Jeffrey C. Lagarias} \\
\vspace*{.2\baselineskip}
AT\&T Labs- Research \\
Florham Park, NJ 07932-0971 \\
\vspace{1\baselineskip}
{\em Yang Wang} \\
School of Mathematics \\
Georgia Institute of Technology \\
Atlanta, GA 30332 \\
\vspace*{1\baselineskip}
\end{center}


%
%
%
%
\section{Introduction}
\hsp

Aperiodic and self-similar structures
in $\RR^d$ have been
extensively studied 
using tilings of $\RR^d$ as models. Among these are
the classes of self-similar tilings
and, more generally,  self-affine tilings.
Such tilings  have been proposed as models for quasicrystalline structures 
(\cite{GahKli97}, \cite{Hof92},  \cite{Jan92}, \cite{LevSte84}, 
\cite{LunPle87}, \cite{Ro96}); they also arise in constructions of 
compactly supported wavelets and multiwavelets (  \cite{FW},\cite{GHR},
\cite{GM92}). 
An alternate
 method of modelling quasicrystalline structures uses discrete sets,
more specifically Delone sets (defined below), see \cite{Lag96},
\cite{LagI}, \cite{LagII}, 
\cite{Mey95}, \cite{Moo95}, \cite{Moo97}, \cite{Sen95}.
These sets model the atomic positions
occupied in the structure. In terms of tilings, such discrete sets  
can be viewed either as
{\em tiling sets},  
representing
the translations used in forming tilings by translation
of a finite number of  different
prototile types, or as {\em control points} marking in some way 
the location of each tile.  

Comparison of these two types of models, which appear
rather different, motivates the question: 
Is there a notion of self-similarity appropriate 
to discrete sets and Delone sets? This paper develops such a notion,
which
is based on a system of functional equations dual to the
functional equations associated to self-affine tilings and multi-tilings.

We first recall the functional equation associated 
to the construction of finite sets of tiles
$\{T_1, \ldots, T_n \}$ which tile $\RR^d$ with special 
kinds of self-affine tilings.
The tiles are solutions to a finite system of set-valued 
functional equations which we call multi-tile equations, which 
encode a self-affine property.
An {\em inflation map} $\phi: \RR^d \to \RR^d$ is an expanding linear map
$\phi ( \bx ) = \he{A} \bx$ in which $\he{A}$ is an expanding $n \times n$
real matrix, i.e. all its eigenvalues $| \lambda | > 1$.
Let
$\{ \sD_{ij} : 1 \le i, j \le n \}$ be finite sets in $\RR^d$ 
called digit sets.

\paragraph{Multi-Tile Functional Equation.}
{\em The family of compact sets
$(T_1, T_2, \ldots, T_n )$ satisfy the system of equations
\begin{equation}\label{101}
\he{A} (T_i) = \bigcup_{j=1}^n (T_j + \sD_{ji} )
\eeq
for $1 \le i \le n$.
}

\vspace*{+.1in}
This functional equation is set-valued, i.e. points are counted without 
multiplicity in the set union on the right side of \eqn{101}.
The {\em subdivision matrix} associated to \eqn{101} is
\beql{102}
\he{S} = [ | \sD_{ij} | ]_{1 \le i,j \le n}.
\eeq

These functional equations have a nice solution theory when the substitution
matrix satisfies the following extra condition.

\begin{defi}
{\em A nonnegative real matrix $\he{S}$ is {\em primitive} if some power 
$\he{S}^k$
has strictly positive entries.}
\end{defi}

In the one-dimensional case the subdivision matrix associated
to \eqn{101} is always primitive. 
It is known that when $\he{S}$ is primitive the functional 
equation \eqn{101} 
has a unique nonempty solution 
$$\sT := (T_1, \ldots, T_n )$$ 
in which all $T_i$
are compact sets
(see \cite[Theorem 2.3]{FW}, \cite{GHR}); in this case
all $T_i$ are nonempty.  In the imprimitive case
there can be more than one  nonempty solution $\sT$ := $(T_1, \ldots, T_n )$
in which all  $T_i$
are  compact sets (see examples in \cite{GHR}),  
and the theory becomes more complicated. There  are however at most 
finitely many such solutions, and under general conditions
there is a unique
``maximal'' compact solution, in which all $T_i$ are nonempty, see \cite{FW}.
In this paper we primarily  consider functional equations where 
the subdivison matrix $\he{S}$ is primitive, though some results
apply for general $\he{S}.$

We are interested
in the case where these sets $T_i$ have positive Lebesgue measure and
can be used in tiling $\RR^n$ by a self-affine tiling. 
The tiling sets for such tilings are special solutions
to a second functional equation involving the same
data, which is ``adjoint'' to the multi-tile
functional equation.  This functional equation counts multiplicities
of sets, unlike \eqn{101}, and we consider solutions to
it that are multisets.

\paragraph{Inflation Functional Equation.}
{\em The multiset family $ \sX := (X_1, X_2, \dots, X_n)$ satisfies 
the system of equations
\beql{103}
X_i = \bigvee_{j=1}^n ( \he{A} (X_j) + \sD_{ij} ) , \quad
1 \le i \le n,
\eeq
where $\sD_{ij}$ are finite sets of vectors in $\RR^d$.
}

Here each $X_i$ is a multiset (as defined in \S2),
 and $\bigvee_{j=1}^n $ denotes
multiset union, as defined in \S2. We can view the inflation
functional equation as a fixed point equation $\psi(\sX)=\sX$,
where $\psi(\cdot)$ is an operator taking multisets to multisets,
defined by $\psi(\sX) :=(X_1^{'}, ... , X_n^{'})$
with 
$$X_i^{'} := \bigvee_{j=1}^n ( \he{A} (X_j) + \sD_{ij} ).$$
The inflation functional equation can also be written
as a system of  equations for the  multiplicity functions:
\beql{104}
m_{X_i} ( \bx ) = \sum_{j=1}^n m_{\he{A} (X_j) + \sD_{ij}} ( \bx )
	= \sum_{j=1}^n \sum_{\bd\in\sD_{ij}}m_{X_j}(\he A^{-1}(\bx -\bd)) , 
\quad
\fa \quad \bx \in \RR^d .
\eeq

A {\em self-replicating multi-tiling} consists of a pair $(\sT, \sX)$ of
solutions to the multi-tile functional equation and
the inflation functional equation such that:

(1) The solution $\sT := (T_1, \ldots, T_n )$ to the multi-tile
functional equation has all sets  
$T_i$ of positive Lebesgue
measure. 

(2) The solution $\sX =(X_1, \ldots, X_n)$ to the inflation
functional equation are sets 
(multi-sets with all multiplicities one)
and $\bigcup_{i = 1}^n (T_i + X_i)$ is a tiling of $\RR^d$,
using the $T_i's$ as prototiles.

\noindent This notion extends the notion of self-replicating tilings studied
in Kenyon~\cite{Ke90}, \cite{Ke92}, which allow only one type
of tile ($n = 1$). Later  studies (\cite{Ke94}, \cite{Ke96})
allowed $n$ tile types but
restricted the inflation matrix $\he{A}$ to be a similarity.

We call a family $\sX$ a {\em self-replicating
multi-tiling family} (for a fixed inflation functional equation)
if it appears as part of a self-replicating multi-tiling $(\sT, \sX)$.
It follows that each multiset
$X_i$ in such a family is a uniformly discrete set.
There are very strong restrictions on the data $(\he{A}, \sD_{ij})$
on the functional equations for a self-replicating multi-tiling
to exist, for example  the Perron eigenvalue condition given below.

This  paper
studies solutions to the inflation functional equation
that are discrete, 
including solutions that do not correspond to tilings.
The inflation functional equation has properties that significantly
differ from  those of the multi-tile functional equation. 
For example, the multi-tile functional equation has a unique
solution (in the case of primitive subdivision matrix)
because the solution $\sT$ is given by the unique fixed point of
a contracting system of mappings. In contrast, the inflation
functional equation involves  an expanding system of mappings, and
its solutions are not compact sets. It may
have infinitely many different solutions,
even infinitely many very nice solutions in some cases. 
Our replacement for the ``contracting'' condition,
is to restrict attention to solutions
$\sX = (X_1, \dots,X_n)$ having the 
property that all
$X_j$ are discrete multisets. 
 A  set ${\us X}$ in $\RR^d$ is {\em discrete} if each bounded set in $\RR^d$
contains finitely many elements of ${\us X}$.
A multiset $X$ is {\em discrete} if its underlying set $\us X$ is
discrete and each element in $X$ has a finite multiplicity.
A multiset family $\sX$ is discrete if each multiset $X_i$ in it is 
discrete.

In \S2  we give precise definitions and statements of
the main results in the paper. Below we summarize
the general contents of the other sections.

In \S3 we  develop a structure theory for multiset solutions to the
inflation functional equation that are discrete. 
These results are proved for general inflation functional
equations, with no primitivity restriction. 
We show that discrete multiset solutions
 decompose uniquely into a finite number of
irreducible discrete multisets, and that each irreducible
discrete multiset is characterized by a finite set of points in it.
However there are significant restrictions on the
inflation functional equation data $(\he{A}, \sD_{ij} )$
necessarry for the existence of any  
discrete multiset solution, as indicated in later sections.

In \S4--\S7  we study multiset solutions 
which correspond more closely to tilings.
A multiset $X$ is {\em weakly uniformly discrete} if there is
a positive radius $r$ and a finite constant $m \geq 1$ 
such that each ball of radius $r$ contains at most $m$ points
of $X$, counting multiplicities; it is 
{\em relatively
dense} if there is a radius $R$ such that each ball of radius
$R$ contains at least one point of $X$.
A multiset $X$ is a {\em weak Delone set} $X$ if it is weakly uniformly
discrete and relatively dense.

A solution $\sX = (X_1, ..., X_n)$
to the inflation functional equation is a {\em weak substitution Delone
multiset family} if each $X_i$ is a weak Delone multiset;
In studying solutions which are weak Delone multisets, 
we restrict attention to
inflation functional equations that have a primitive
subdivision matrix $\he{S}$. Then we make use of
Perron-Frobenius theory, which
asserts that a primitive  nonnegative real matrix $\he{M}$
has a positive real eigenvalue $\lambda ( \he{M} )$ 
such that:
\begin{itemize}
\item[(i)]
$\lambda (\he{M})$ has multiplicity one.
\item[(ii)]
$ \lambda ( \he{M} ) > | \la ' |$ for all eigenvalues $\la'$ of $\he{M}$ 
with $\la' \neq \la(\he{M})$.
\item[(iii)]
${\he M}$ has  both right and left eigenvectors 
for eigenvalue  $\lambda ( \he{M} )$ 
which have  positive real entries.
\end{itemize}
We call $\lambda ( \he{M} )$ the {\em Perron eigenvalue} of $\he{M}$;
it is equal to the spectral radius of $\he{M}$.

For inflation functional equations whose substitution
matrix is primitive, we 
show that those multiset equations having a weak Delone set
solution must satisfy the following: 

\noindent{\bf Perron eigenvalue condition}: 
{\em The Perron eigenvalue $\lambda(\he{S})$ of
the subdivision matrix $\he{S}$ satisfies}
\beql{eq105}
\lambda(\he{S}) = |\rdet (\he{A})|.
\eeq

More generally, there are inequalities relating the Perron
eigenvalue and properties of
solutions of the multi-tile and inflation functional equations,
as follows:

(1) A necessary condition for the multi-tile functional
equation with primitive subdivision matrix
to have $\sT = (T_1, ..., T_n)$ with
some (and hence all) $T_i$ of positive Lebesgue measure is that
$$ \lambda(\he{S}) \geq |\rdet (\he{A})|. $$

(2) A necessary condition for the inflation functional
equation with primitive subdivision matrix to have
a solution $\sX = (X_1, ..., X_n)$ with some (and hence all) $X_i$ 
weakly uniformly discrete is that 
$$ \lambda(\he{S}) \leq |\rdet (\he{A})|. $$

Inequality (1) is established by taking the Lebesgue measure on both 
sides of the multi-tile
equation (\ref{101}), see Theorem~\ref{th51a}.
 In the  special case when $\he{A}$ is a similarity
Mauldin and Williams~\cite{MW88} showed the stronger result that if  
$\lambda(\he{S}) < |\rdet (\he{A})|$ then 
the Hausdorff dimension of each $T_i$ must be less than $d$. 
Inequality (2) is established in Theorem~\ref{thm35}.

In \S5, under the assumption that the Perron
eigenvalue condition holds, we give some
necessary and sufficient conditions for the associated
multi-tiling equation to have solutions $T_i$  having
positive Lebesgue measure. Then in \S6 we use these to show
that this condition is equivalent to 
existence of
weak substitution Delone set solutions 
$\psi(\sX) = \sX)$ is closely associated with
self-replicating multi-tilings. 
In particular, the  associated multi-tile functional
equation necessarily has a compact solution $\sT= (T_1, T_2, \dots, T_n)$
with the $T_i$ all having positive Lebesgue measure,
see Theorem~\ref{th14}, and some iterate $\psi^N(\cdot)$
of the inflation functional equation then has a solution that is a
self-replicating multi-tiling. In \S7 we supplement this
with a sufficient condition to have a self-replicating
multi-tiling.

In \S8 we give examples showing the limits of our results.

We distinguish  {\em multi-tilings}, which are tilings of
$\RR^d$ using several types of tiles, from
{\em multiple tilings}, which are arrangements of tiles in
$\RR^d$ such that almost all points in $\RR^d$ are
covered exactly $p$ times, for some positive integer $p$. We
call $p$ the {\em thickness} of the multiple tiling, and
sometimes call such a tiling {\em $p$-thick.} A $1$-thick tiling
is just an ordinary tiling. The constructions
of this paper can have associated to them tilings which are multiple in
both senses, that is,  $p$-thick multi-tilings for some $p \ge 2$. 
By definition self-replicating multi-tilings are $1$-thick
tilings. We remark that the tiling set of a $p$-thick
tiling cannot in general be partitioned to give a union of $1$-thick tilings,
even in dimension $1$.

To conclude this introduction, given
 any Delone set $X$ one can associate a topological dynamical system
$([[X]], \RR^d)$ with an $\RR^d$-action,
in which $[[X]]$ is the closure of the orbit of $X$ 
under the $\RR^d$-translation action
in an appropriate topology, see Solomyak~\cite{Sol98b}. 
For substitution Delone set families these dynamical systems can 
be viewed as a 
generalization of substitution dynamical systems
(Queffelec~\cite{Que87}), in that every 
primitive substitution dynamical system is topologically conjugate
to a suitable substitution Delone set dynamical system.
Under sufficiently strong extra
hypotheses  substitution Delone set dynamical systems
are minimal and uniquely ergodic. We hope to return to this
question elsewhere.

\paragraph{Acknowledgments.} We are indebted to the reviewers
for very detailed comments improving the paper.

\paragraph{Notation.} 
The positive integer $d$ refers to the dimension $\RR^d$.
We use $B_R(\bx) := \{ \by \in \RR^d: || \by - \bx|| \le R \}$
to denote the closed Euclidean ball of radius $R$ around $\bx$.
The positive integer $n$ denotes the number of tile types or terms in
the inflation functional equation. The positive integer $p$ refers to the
period of a cycle, and also sometimes to the 
the thickness of a 
multiple covering or tiling  of $\RR^d$.
%
%
%
%
\setcounter{section}{1}
\setcounter{equation}{0}
\section{Definitions and Main  Results}

We consider solutions to the inflation functional equation that
are multisets.

\begin{defi}
{\rm \begin{itemize}
\item[(i)] A {\em multiset} $X$ in $\RR^d$ is 
a set ${\us X}$ together with a positive integer valued function,
the multiplicity function 
$m_X: {\us X} \to \ZZ_{> 0}$, 
which assigns to each element $\bx \in {\us X}$ 
a ``multiplicity'' $m_X(\bx)$. 
Given a multiset $X$ we shall use
$\us X$ to denote the underlying set of $X$, i.e. without counting the
multiplicity,  and use either  $m_X(\bx)$ or $m(X, \bx)$ 
to denote the multiplicity function.
A multiset is called an {\em ordinary set} 
if every element in $X$ has multiplicity $1$. We extend the
multiplicity function to all $\bx \in \RR^d$ by setting
$m_X(\bx) = 0$ if $\bx \not\in {\us X}$.

\item[(ii)] Given multisets $X_1$ and $X_2$ in $\RR^d$ 
we say $X_1 \subseteq X_2$
provided $\us X_1 \subseteq \us X_2$ and $m_{X_1}(\bx) \leq m_{X_2}(\bx)$
for all $\bx \in X_2.$ In particular,  ${\us X} \subseteq X.$
\end{itemize}
}
\end{defi}

It is 
convenient to extend the definition of multiplicity function
$m_X(\bx)$ to all $\bx \in \RR^d$, by setting $m_X(\bx)=0$
if $\bx \not\in {\us X}$.

The term ``multiset'' is attributed  to
N. G. deBruijn
(see Knuth \cite[p. 636]{Kn81}) and the concept has many uses.
The combinatorial view of
a multiset is as a collection of possibly repeated elements.
For example, $X= \{0,0,1,3,4,4,4\}$ represents a multiset in $\RR$
in which $0$ is counted
twice and $4$ is counted 3 times. Thus $\us X = \{0,1,3,4\}$ and
$$
    m_X(0)=2,~~ m_X(1)=m_X(3) = 1,~~m_X(4)=3.
$$

One may
 alternatively regard a multiset $X$ in $\RR^d$ as a nonnegative
integer-valued pure discrete measure on $\RR^d$. 
The results of this paper could be formulated in
 a measure-theoretic framework,
and the  inflation functional equation then expresses an
equality of measures. 
A reviewer observed that a measure-theoretic treatment might
prove useful for further generalizations, 
but we do not attempt it here.

\begin{defi}
{\em
For any multisets $X$ and $Y$ the {\em multiset union} $X \vee Y$ 
is the multiset having
the multiplicity function
\beql{101b}
m_{X\vee Y}  = m_X +m_Y,
\eeq
and the
{\em multiset intersection} $X \wedge Y$ is the
multiset having  multiplicity function
\beql{101c}
    m_{X\wedge Y} = \min\{m_X, m_Y\}.
\eeq
For a multiset $X$ and a set (or multiset) $\sD$ the {\em multiset sum} 
$X+\sD$ is
\beql{101d}
X+\sD := \bigvee_{\bd\in\sD} (X+\bd). 
\eeq
}
\end{defi}

\begin{defi}
{\em
A {\em multiset family}  $\sX=(X_1, \dots, X_n)$ is a finite
vector of multisets $X_i$.
We call a multiset family $\sX=(X_1, \dots, X_n)$ an {\em $n$-multiset family}.
For $n$-multiset families $\sX^{(1)}$ and $\sX^{(2)}$ we define
$$
	\sX^{(1)} \vee \sX^{(2)} = (X_1^{(1)}\vee X_1^{(2)}, \dots,
	X_n^{(1)}\vee X_n^{(2)}).
$$
}
\end{defi}

\begin{defi}\label{de21b}
{\em
A multiset family $\sX= (X_1, ..., X_n)$ is {\em discrete} if for
each $i$ the multiset $X_i$ is discrete, i.e. the underlying set 
$\underline{X}_i$ is discrete and elements
in $X_i$ have finite multiplicity. 
}
\end{defi}

In \S3 we develop a structure theory
for discrete multiset families, which decomposes them
into irreducible families. Given
an inflation functional equation,
let  $\psi(\cdot)$ be the inflation operator associated to it,
which takes  $n$-multiset families to $n$-multiset families, 
as defined in \S3.

\begin{defi}\label{de21a}
{\em 
(i) A multiset family $\sX= (X_1, ..., X_n)$ 
satisfying an inflation functional equation
$\psi ( \sX ) = \sX$ is {\em indecomposable} if it cannot be
partitioned as $\sX = \sX^{(1)} \vee \sX^{(2)}$ with each $\sX^{(i)}$
nonempty such that
$\psi ( \sX^{(i)}) = \sX^{(i)}$ for each $i=1,2$. Otherwise
it is decomposable.

(ii) A multiset family $\sX= (X_1, ..., X_n)$ 
satisfying an inflation functional equation
$\psi ( \sX ) = \sX$ is {\em irreducible} if there does
not exist a nonempty $\sX'= (X_1^{'}, ..., X_n^{'})$ 
satisfying $\psi ( \sX^{'}) = \sX^{'}$
with  $\sX^{'} \subsetneq \sX$ in the sense that  $X_i^{'} \subseteq X_i$
for all $i$, with some $X_i^{'} \neq X_i$. Otherwise it is
reducible.
}
\end{defi}

The  concepts of indecomposable and irreducible 
multiset family (for a fixed inflation functional equation) are equivalent.
If $\sX$ is  irreducible then it is necessarily indecomposable.
Conversely, if $\sX$ is reducible, with $\sX' \subsetneq \sX$,
and $\psi (\sX') = \sX'$, 
then $\sX''= \sX - \sX'$ (subtraction done on multiplicity
functions) is a multiset family with  $\sX = \sX' \bigvee \sX''$, and 
linearity of the inflation functional
equation yields $\psi(\sX'')= \sX''$, so that 
$\sX$ is decomposable.

\begin{theorem}[Decomposition Theorem]\label{th11}
Let $\sX$ be a multiset family $\sX$ that is discrete and 
satisfies an inflation functional equation $\psi ( \sX ) = \sX$. 
Then $\sX$ uniquely partitions into a finite number of irreducible
discrete multiset families that satisfy the same inflation functional equation.
\end{theorem}

An inflational functional equation $\psi ( \sX ) = \sX$ may not have 
a discrete solution, see Example \ref{no-discrete}.
In \S3 we characterize indecomposable discrete multiset
families $\sX$ satisfying inflation functional equations,
as follows.

\begin{theorem}~\label{th12}
Let $\sX$ be a multiset family which satisfies
an inflation functional equation $\psi(\sX) = \sX$
and is discrete and indecomposable. Then $\sX$  is 
irreducible and is  generated
by a finite ``seed'' $\sS^{(0)}= (S_1, S_2, ..., S_n)$, which
consists of a periodic cycle $Y = \{ (\bx_j, i_j) : 1 \le j \le p \}$
in which $\bx_j \in X_{i_j}$ and there is some 
$\bd_j \in  \sD_{i_{j+1}, i_{j}}$
with 
$\bx_{j + 1} = \he{A} \bx_j + \bd_j$, with
$ (\bx_{p+1}, i_{p+1}) = (\bx_1, i_1).$ That is, 
$$\sX = \lim_{N \to \infty} \phi^N(\sS^{(0)}). $$
The periodic cycle $Y$ is the only periodic cycle in $\sX$ and its
elements each have multiplicity one.
\end{theorem}

This result appears,
in a more precise form, as Theorem~\ref{th21}.
For each $p$ there are only finitely many periodic cycles $Y$
and they can be effectively enumerated.
However not  all periodic cycles $Y$ generate irreducible discrete
multiset families. 
We show that there  is an
algorithmic procedure, which when given  any such ``seed'' as input,
has one of three outcomes:

(1) If the generated multiset system
 $\sX$ is discrete and irreducible, the procedure  eventually
halts and certifies this holds.

(2) If the generated multiset system $\sX$ is not irreducible,
or if the limit does not exist,
the procedure  eventually halts and certifies this holds.

(3) If the generated multiset system $\sX$ is irreducible and
not discrete, the procedure may not halt.

We also prove in \S3 a dichotomy concerning
 the multiplicities of elements
appearing in the multisets $X_i$ in an irreducible discrete
multiset family $\sX$ satisfying an inflation functional equation
having a  primitive subdivision matrix :
Either they are all bounded by the period $p$ of the
generating cycle or else they all have unbounded multiplicities.
(Theorem~\ref{th22}.) We do not know of
 an effective computational
procedure to tell in general which of these alternatives occurs.

 In \S4-\S7  we restrict to the case of inflation functional
equations having primitive subdivision matrix.
We study Delone
set solutions and self-replicating multi-tilings.

\begin{defi}~\label{de108b}
{\em 
(i) A multiset $X$ is {\em weakly uniformly discrete}
if there is a radius  $r > 0$ and a finite $m \ge 1$
such that any open ball of radius $r$ contains at 
most $m$ points of $X$, counting points with multiplicity.

(ii) A multiset $X$ is {\em uniformly discrete}
if there is a radius  $r > 0$ such that 
such that any open ball of radius $r$ contains at 
most one point of $X$, counted with multiplicity.
Such an $X$ is necessarily an ordinary set.

(iii) A multiset $X$ is {\em relatively dense} if
there is a finite radius $R >0$ such that each
closed ball of radius $R$ contains at least one
point of $X$.
}
\end{defi}

Note that  a  finite union of weakly uniformly discrete sets
is weakly uniformly discrete. 

\begin{defi}~\label{de109}
{\em
(i) A multiset $X$ is a {\em weak Delone multiset} 
if it is weakly uniformly discrete and relatively dense.

(ii) A multiset $X$ is a {\em  Delone set} 
if it is uniformly discrete and relatively dense.
That is, it is an ordinary set which has both a finite nonzero packing
radius and a finite covering radius  by equal spheres.
}
\end{defi}

A finite union of Delone sets need not be a Delone set, while
a finite union of weak Delone multisets is a weak Delone multiset.

\begin{defi}\label{de111}
{\em 
(i) A multiset family 
$\sX= (X_1, \dots, X_n)$ is 
a {\em weak Delone multiset family} 
if $\sX$
satisfies the inflation functional equation
\eqn{103} in $\RR^d$ and
each multiset
$X_i$ is a weak Delone multiset.

(ii) It is a {\em substitution Delone set family} if it is
a weak substitution Delone multiset family and 
each $X_i$ is a Delone set.
}
\end{defi}

Weak substitution Delone multiset families can sometimes
be associated to $p$-thick multiple tilings, for some $p \ge 2$.

In \S4 we prove the following result.

\begin{theorem}[Perron Eigenvalue Condition]\label{th13}
If the inflation functional equation $\psi(\sX) = \sX$
with primitive subdivision matrix $\he{S}= [|\sD_{ij}|]$
has a solution $\sX$ that is a weak Delone multiset family, then the
Perron eigenvalue $\lambda(\he{S})$ of $\he{S}$ satisfies
$$ \lambda(\he{S})= |\rdet (\he{A})|. $$
\end{theorem}

The narrowest  class of solutions to the inflation functional
equation are those having the following tiling property.

\begin{defi}\label{de12}
 {\em
A multiset family $\sX =  (X_1, \dots, X_n)$ is called
a {\em self-replicating multi-tiling family} for a given system
$(\he{A}, \sD_{ij} )$ in \eqn{103} with primitive subdivision 
matrix $\he{S}$ if
\begin{itemize}
\item[\rm (i)]
$\sX$ is a substitution Delone set family for $( \he{A}, \sD_{ij} )$.
\item[\rm (ii)]
The associated multi-tile equation \eqn{101} has a unique solution 
$\sT :=(T_1, \ldots, T_n )$ with each $T_i$ of positive Lebesgue measure.
\item[\rm (iii)]
The sets $\{ T_i + \bx_i : 1 \le i \le n \quad {\rm and} \quad
\bx_i \in X_i \}$ tile $\RR^d$.
\end{itemize}
}
\end{defi}

This  definition requires that the multiset family $\sX$ give
a $1$-thick tiling, rather than a $p$-thick
multiple tiling for some $p \geq 2$, or
a $p$-packing for some $p \geq 1$ that is not a tiling. 

In \S5 we study the multi-tiling functional equation and
show that a necessary
condition for a solution with some $T_i$ having positive
Lebesgue measure is that  $\lambda(\he{S}) \ge |\det(\he{A})|$.
Assuming that the Perron eigenvalue
condition  $\lambda(\he{S}) = |\det(\he{A})|$ holds,
in    Theorem~\ref{th44} we characterize when  solutions 
of positive Lebesgue measure exist.

In \S6 we prove the following result, which shows that  substitution
Delone set families  are related to the existence of self-replicating
multi-tilings.
In this result  $\psi^N(\cdot)$
denotes the inflation operator  $\psi(\cdot)$ composed with
itself $N$ times, whose associated inflation matrix is $\he{A}^N$.

\begin{theorem}\label{th14}
Let $\psi ( \sX ) = \sX$ be an inflation functional equation 
that has a primitive subdivision matrix.
Then the following conditions are equivalent:
\begin{itemize}
\item[\rm (i)]
For some $N >0$ there exists a weak Delone multiset family $\hat{\sX}$
satsfying $\psi^N (\hat \sX ) = \hat\sX$.
\item[\rm (ii)]
For some $N > 0$ there exists a self-replicating multi-tiling family
$\hat\sX$ for the inflation functional
equation   $\psi^N(\cdot)$. 
\item[\rm (iii)] 
The subdivision matrix $\he{S}$  satisfies the Perron eigenvalue condition 
$\lambda(\he{S}) = |\rdet(\he{A})|$ and 
the unique compact solution $(T_1, \ldots, T_n )$ 
of the associated multi-tile functional equation 
consists of sets $T_i$ that have positive Lebesgue measure, 
$1 \le i \le n$, and each $T_i$ is the closure of its
interior $T_i= \overline{T_i^{\circ}}.$
\end{itemize}
\end{theorem}

These conditions imply that the Perron eigenvalue condition holds.
If we instead assume that the Perron eigenvalue condition holds
then condition (i) can be relaxed to assuming only the existence 
of a weakly uniformly discrete solution, 
as shown in  Theorem~\ref{th41}.
It is known there are many restrictions on the data
$ ( \he{A} , \sD_{ij} )$ in order for  the conditions
of Theorem~\ref{th14} to hold.
A number of different necessary (or sufficient) conditions on 
the data $ ( \he{A} , \sD_{ij} )$ are given
in  Kenyon \cite{Ke90}-\cite{Ke96},
Lagarias and Wang \cite{LW96a}- \cite{LW97},
Praggastis \cite{Prag}, \cite{Pr99} and 
Solomyak \cite{Sol97}--\cite{Sol98b}.

In \S7 we give a sufficient condition for a substitution
Delone set family to be a self-replicating multi-tiling 
(Theorem~\ref{th51}). 
In \S8 we give examples and counterexamples showing the limits 
of our results.

%
%
%
%
\setcounter{section}{2}
\setcounter{equation}{0}

\section{Structure of Discrete Multiset Solutions}

\hsp
We consider multiset families $\sX = (X_1, \ldots, X_n )$ 
satisfying the inflation functional equation
\beql{201}
X_i = \bigvee_{j=1}^n ( \he{A} (X_j) + \sD_{ij} ), \quad 1 \le i \le n.
\eeq
The
 {\em inflation operator} $\psi(\cdot)$ associated
to data $(\he{A}, \sD_{ij})$ is an operator which 
maps arbitrary  $n$-multiset families $\sX$ to $n$-multiset families
$\psi(\sX) = \sX^{'}$, as follows.
Given $\sX$ (not necessarily satisfying (\ref{201})), we 
define  $\psi(\sX) = ( X_1^{'} , \ldots, X_n^{'})$ by
\beql{202}
X_i^{'} :=
\bigvee_{j=1}^n ( \he{A} (X_j) + \sD_{ij}),  \quad 1 \le i \le n.
\eeq
The inflation functional equation \eqn{201} then asserts that
$\sX$ is a fixed point of $\psi ( \cdot )$,
i.e. $\psi ( \sX )=\sX$.

In this section we determine the structure of 
multiset families $\sX$ that are
discrete and  satisfy $\psi ( \sX ) = \sX$, where we  put
no restriction on
$\psi(\cdot)$. The results allow imprimitive substitution matrices,
unless otherwise stated.
We also describe the structure of such sets that are also 
irreducible
(Theorems~\ref{th21} and \ref{th22}). We  show that every
multiset family $\sX$ that is discrete and
satisfies $\psi ( \sX ) = \sX$ uniquely
partitions into a finite number of irreducible sets (Theorem~\ref{th11}).

We consider solutions to the inflation functional equation \eqn{201}
built up by an iterative process starting from a finite multiset
which is a ``seed.''
Let
\beql{204}
\sS^{(0)} := ( {S}_1^{(0)} , \ldots,  S_n^{(0)} )
\eeq
be a system of finite multisets, and iteratively define the 
finite multiset family
$\sS^{(k+1)}$ by $\sS^{(k+1)} =\psi ( \sS^{(k)} )$, i.e.
$$
\sS_i^{(k+1)} := \bigvee_{j=1}^n ( \he{A} (S_j^{(k)}) + \sD_{ij}),  
\quad 1 \le i \le n.
$$
We say that $\sS^{(0)}$ satisfies the {\em inclusion property} 
$\sS^{(0)} \subseteq \sS^{(1)}$ if
\beql{205}
 S_i^{(0)} \subseteq  S_i^{(1)} \quad {\rm for} \quad 1 \le i \le n ~,
\eeq
in the sense of multisets.

\begin{lemma}\label{le21}
If a multiset family $\sS^{(0)}$ is finite and 
satisfies the inclusion condition
$\sS^{(0)} \subseteq \psi ( \sS^{(0)} )$ then
\beql{206}
\sS^{(k)} \subseteq \sS^{(k+1)} \quad \mbox{for all}\quad k \ge 0 ~.
\eeq
The following limit set is well-defined
\beql{207}
\us X_i := \lim_{k \to \infty} \us {S}_{i}^{(k)} , \quad 1 \le i \le n ~,
\eeq
where $\us X_i$ is regarded as a countable point set. For each
$\bx \in {\us X_i}$ the multiplicity function has a limit,
\beql{207a}
m_{X_i}( \bx) :=  \lim_{k \to \infty} m({S}_i^{(k)},\bx), 
\eeq
where this multiplicity may take the value $+\infty.$
If all multiplicities in (\ref{207a}) remain finite, then the
multisets $X_i$ are well-defined and the
multiset family $\sX = ( X_1 , \ldots, X_n )$ satisfies
$\psi ( \sX ) = \sX$.
\end{lemma}

\paragraph{Remark.}
If all the multiplicities remain finite, we write 
$\sX := \lim_{k \to \infty} \sS^{(k)}$.
The limit multiset family $\sX$ need not be discrete, and the underlying
sets could even be dense.
\paragraph{Proof.}
The inclusion \eqn{206}
follows by induction on $k$, since $\sS^{(k-1)} \subseteq \sS^{(k)}$ yields
$$
\sS^{(k)} = \psi ( \sS^{(k-1)} ) \subseteq 
\psi ( \sS^{(k)} ) = \sS^{(k+1)} ~.
$$
The other statements follow easily from \eqn{206}.
\eproof

To further analyze multiset families $\sX$ that satisfy
$\psi ( \sX ) = \sX$, we put the structure of a colored directed 
graph on $\sX$. 
The vertices of this graph are the points in the disjoint union of the 
underlying sets
$\us X_i$, with vertices being labelled $(\bx, i)$,
where $\bx \in \us X_i \in \RR^d$ and $i$ is the ``color''.
For each $\bx_j \in X_j$ and each $\bd_j \in \sD_{ij}$, set
$$
	\bx_i' := \he{A} \bx_j + \bd_j \in X_i ~,
$$
and put a directed edge $\bx_j \to \bx_i'$ in the graph. We call
$\bx_i'$ an {\em offspring} of $\bx_j$, and  call $\bx_j$ a
{\em preimage} of $\bx_i'$.
We denote this (infinite) colored  directed graph by $\sG ( \sX )$. 
We also assign to each vertex $\bx \in \us X_i$ a weight which is its
multiplicity $m_{X_i}(\bx)$.

\begin{lemma}\label{le22}
Let $\sX=(X_1,..., X_n)$ be a multiset family 
which is discrete and satisfies the 
inflation functional equation $\psi ( \sX ) = \sX$.
\begin{itemize}
\item[\rm (i)]
There is a finite multiset family $\sS^{(0)} \subseteq \sX$ 
with the inclusion property, such that
\beql{208}
\sX = \lim_{k \to \infty} \sS^{(k)} ~.
\eeq
\item[\rm (ii)]
The directed graph $\sG(\sX)$ contains a directed cycle.
\end{itemize}
\end{lemma}
\Proof
(i)   We first show that there exists a radius 
$R' >0$, depending only on $\psi(\cdot)$,
such that every vertex in $\sX$ can be reached by a directed path in the 
graph $\sG ( \sX )$ from some point $\bx \in \sX \cap B_{R'}(\bo),$
where $B_{R'}(\bo) := \{ \bx \in \RR^d: ||\bx|| \leq R'\}.$  
To show this, we set
\beql{209}
\la ( \he{A}^{-1} ) = \min \{ | \la |~:~\la ~\mbox{ is an eigenvalue of 
$\he{A}^{-1}$} \},
\eeq
and $\la ( \he{A}^{-1} ) < 1$ since $\he{A}$ is expanding.
Fix a $ \rho >1$ such that $\la( \he{A}^{-1} )<1/\rho<1$. It
is shown in Lind~\cite{Lin82} that 
there exists a norm $\| \cdot \|_{\he{A}}$ on $\RR^d$ with the property that
\beql{210}
\| A^{-1} \bx \|_{\he{A}} \le \frac{1}{\rho} \| \bx \|_{\he{A}} 
\quad\mbox{for all} \quad \bx \in \RR^d .
\eeq
The norm is defined by 
$\| \bx \|_{\he{A}} := \sum_{k=1}^\infty \rho^k \|A^{-k} \bx \|$.

    We claim that there exists an $R>0$ such that for any vertex
$\bx_i' \in \us X_i$ with $\|\bx_i'\|_{\A}\geq R$, any
preimage $\bx_j \in \us X_j$ of $\bx_i'$ must satisfy
$$\|\bx_j\|_{\A}  \leq \frac{2}{1+\rho}\|\bx_i'\|_{\A}.$$
To see this, let 
\beql{214}
    C = \max \{\| \bd_j \|_{\he{A}} :~\bd_j \in \sD_{ij}, ~ 1 \le i,j \le n \}
\mbox{~~and~~} R:= \frac{\rho+1}{\rho-1} C.
\eeq
It follows from
$ \bx_j = \he{A}^{-1} ( \bx_i' - \bd_j )$ that
\beql{216} 
\|\bx_j \|_{\he{A}} \le \frac{1}{\rho} 
\left( \| \bx_i' \|_{\he{A}} + C \right) 
 \leq \frac{2}{1+\rho} \| \bx_i' \|_{\he{A}}.
\eeq

Thus if a vertex $\bx$ has $\| \bx \|_{\he{A}} \ge R$ 
then every preimage of $\bx$ in the 
graph $\sG ( \sX )$ is smaller in the norm
$\| \cdot \|_{\he{A}}$ by a multiplicative constant 
$\frac{2}{1+ \rho} <1$.
It follows that every vertex $\bx$ in $\sG(\sX)$ can be reached by a finite 
directed
path in $\sG ( \sX )$ starting from some vertex  $\bx'$ with
$\| \bx' \|_{\he{A}} \le R$. Since the set
$$\tilde{B}_R :=\{\by \in \RR^d : ||\by||_{\he A} \le R\},$$ 
is compact, it is contained in some ball of radius $R'$ around ${\bf 0}$, 
which gives the result.

We now take the multiset family
$\sS^{(0)} = (S_1^{(0)}, \ldots, S_n^{(0)})$ such that
$S_i$ consists of those elements $\bx$ of $X_i$ with
$\|\bx\|_{\A} \leq R$,  
counting  multiplicities, i.e. we take
$$
m_{S_i^{(0)}}(\bx) := m_{X_i}(\bx)~~~\mbox{if}~~ \bx \in {\us S_i^{(0)}}.
$$
The multiset family  $\sS^{(0)}$ is a finite family since $\sX$
is discrete. The inclusion property
\beql{218}
\sS^{(0)} \subseteq \sS^{(1)} = \psi ( \sS^{(0)} )
\eeq
holds because all preimages of vertices of $\sG(\sX)$ in
$\{ \bx \in \RR^d : \| \bx \|_{\he{A}} \le R \}$ lie in this set,
since
$\| \bx_j \|_{\he{A}} \le \frac{1}{\rho} (R+C) \leq R$.
Since
$\sS^{(0)} \subseteq \sX$ and
$\psi ( \sX ) = \sX$ we obtain
$\sS^{(k)} \subseteq \sX$ for all $k$,
viewed as multisets, so that
$$\sX ' = \lim_{k \to \infty} \sS^{(k)}$$
exists and is a multiset.
Now $\sS^{(0)}$ contains all points $\bx$ of  $\sX$ with the correct
multiplicities for $\| \bx \|_{\he{A}} \le R$.
Now, by induction on $k \ge 0$, one proves:

For all $\bx$ with 
$\| \bx \|_{\he{A}} \leq \left( \frac{1+ \rho}{2} \right)^k R$, one has
for $1 \le i \le n$ that the multiplicity
$$m(\sS^{(k)}_i,\bx) = m_{X_i}(\bx).$$

The base case $k=0$ holds by construction, and the induction
step follows because all preimages of a vertex in 
$\{\|\bx \| \leq  (\frac{1+\rho}{2})^k R \}$ 
lie in $\{ \| \bx \| \leq (\frac{1+ \rho}{2})^{k-1} R \}$.
We conclude that $\sX ' = \sX$, which gives (i).

(ii) We choose $R$ as in (i). Pick a point $\bx$ in $\sX$
with 
$$\bx \in \tilde{B}_R := \{\by \in \RR^d : ||\by||_{\he A} \le R\},$$ 
which must exist by the
argument in (i). All preimages of $\bx$ necessarily lie in $\tilde{B}_R$.
Since
$\sX$ is discrete, if
 we follow a chain of successive preimages of a vertex $\bx$ in 
$\sG(\sX)$
 we stay in the  finite set  $\sX \cap \tilde{B}_R$, so some vertex
must occur twice. The path from this vertex to itself forms
a directed cycle in $\sG ( \sX )$.
\eproof

Recall that a multiset family $\sX$ satisfying
an inflation functional equation is {\em indecomposable} if it
has no nontivial partition $\sX = \sX_1 \vee \sX_2$, with
both $\sX_i$ satisfying the same functional equation,
and it is {\em irreducible} if there is no nonempty
$\sX' \subsetneq \sX$ satisfying the same functional equation.

\begin{theorem}[Irreducible Set Characterization]\label{th21}
Suppose that the multiset family $\sX$ is discrete and satisfies an inflation
functional equation $\psi ( \sX ) = \sX$.
The following are equivalent:
\begin{itemize}
\item[\rm (i)] $\sX$ is indecomposable.
\item[\rm (ii)] $\sX$ is irreducible.
\item[\rm (iii)]
    The graph $\sG ( \sX )$ contains exactly one directed cycle (which 
may be a loop), and the elements of this cycle have multiplicity one.
\end{itemize}
If condition {\rm (iii)} holds, let $Y=\{\bx_1, \dots, \bx_p\}$ 
be the cycle in $\sG(\sX)$
with multiplicity one and define the multiset family 
$\sS^{(0)} = (S_1^{(0)}, \ldots, S_n^{(0)} )$ by 
$S_i = \{\bx\in Y:~ \bx \mbox{~has color~} i\}$. 
Then $\sS^{(0)} \subseteq \psi ( \sS^{(0)} )$ and
\beql{219}
\sX = \lim_{k \to \infty} \sS^{(k)}.
\eeq
\end{theorem}
\paragraph{Remark.}
In view of (iii) and \eqn{219} we call $\sS^{(0)}$ the {\em generating cycle} 
of the irreducible
multiset system $\sX$, and we call $p$ the {\em period} of this cycle.

\paragraph{Proof.}
We show  (i) $\Leftrightarrow$ (ii)  $\Rightarrow$ (iii) $\Rightarrow$ (i).
Just before stating Theorem~\ref{th11}
we showed that  implication (i) $\Leftrightarrow$ (ii) holds for
all multiset solutions to a fixed inflation functional equation.

(ii) $\Rightarrow$ (iii).
Lemma~\ref{le22} shows that $\sX$ is generated by the points
$\bx$ in $\sX$ that lie in the compact region $\| \bx \|_{\he{A}} \le R$,
and that $\sG(\sX)$ contains vertices forming
 a directed cycle $Y$ inside this region.
The vertices of $Y$ 
define a finite multiset family $\sS^{(0)} = (S_1^{(0)}, \ldots, S_n^{(0)})$,
with  $ S_i^{(0)} = \us X_i \cap Y$, and all elements of $ S_i^{(0)}$
have multiplicity one.
It is clear that $\sS^{(0)} \subseteq \sS^{(1)} = \psi ( \sS^{(0)} )$ 
since each element of $\sS^{(0)}$ has a preimage in $\sS^{(0)}$.
Since each $\sS^{(k)} \subseteq \sX$ the 
limit
\beql{203aa}
	\sX ' := \lim_{k \to \infty} \sS^{(k)}
\eeq
exists and is a multiset family with
$\sX ' \subseteq \sX$.
We have $\sX ' = \sX$, for
if not, this contradicts the irreducibility of $\sX$.

We next show that $\sX$ has multiplicity one on the vertices of $Y$.
Suppose not. The elements of $\sS^{(0)}$ have multiplicity one, and
at some later stage some vertex  $(\bx, i)$ of the cycle $Y$ has 
multiplicity exceeding one in some  $\sS^{(k)}$. Thus there exists a 
path to $(\bx, i)$ of length $k$ arising from some vertex $(\by, j)$
in  $\sS^{(0)}$, which does not arise purely from moving on the
cycle $Y$. Since the vertex $(\by, j)$ can be reached from  $(\bx, i)$
by moving around the cycle $Y$, taking $l$ steps, say, we obtain a
directed path from $(\bx, i)$ to itself of length $k+l$ which does
not stay in the cycle $Y$. Now we obtain two distinct directed paths of
length $(k+l)p$ from  $(\bx, i)$ to itself, one by wrapping around
the cycle $Y$ $k+l$ times, the other by repeating the path of
length $k+l$ $p$ times. Since we can
concatenate these two distinct paths in any
order, it  follows that the multiplicity of
  $(\bx, i)$  in $\sS^{((k + l)pm)}$ is at least $2^m$. This 
implies that the vertex 
 $(\bx, i)$ has unbounded multiplicity as $m \to \infty,$
so that $\lim_{k \to \infty} \sS^{(k)}$
does not exist, a contradiction. Thus $\sX$ has multiplicity one on
the cycle $Y$.

Finally we show that $\sX$ contains only one directed cycle. Suppose not,
and that $\sG ( \sX )$ contains a directed cycle $Y'$ different
from  $Y$. The same argument as above says that $\sX$ is generated
by the cycle $Y'$ and that all elements of $Y'$ in $\sX$ have
multiplicity one. By exchanging $Y$ and $Y'$ if necessary we
may suppose that $Y$ has a vertex not contained in $Y'$.
Now following a path from a vertex in $Y\setminus Y'$ to $Y'$ we see that
there exists an $\bx'\in Y'$ that has a preimage not in $Y'$.
But $\bx'$ also has a preimage in $Y'$ because $Y'$ is a cycle.
This implies that the multiplicity of $\bx'$ is at least 2, a contradiction.
Thus $\sG ( \sX )$ contains exactly one directed cycle.

\noindent
(iii)~$\Rightarrow$~(i).
Suppose that $\sG ( \sX )$ contains exactly one directed cycle $Y$,
of multiplicity one.  We argue by contradiction. Suppose that
$\sX$ were not indecomposable.
Then we can write  
$\sX = \sX ' \vee \sX''$, where both $\sX'$ and  $\sX''$
satisfy the inflation functional
equation.  Now  Lemma~\ref{le22} (ii) implies that 
both $\sG(\sX')$ and 
$\sG ( \sX '' )$ contain a directed cycle.
Therefore $\sG ( \sX )$ either contains two directed cycles,
or else contains a directed cycle of multiplicity at least two.
This contradicts the hypothesis. 
\eproof

 We can now prove the Decomposition  Theorem.

\noindent\paragraph{Proof of Decomposition Theorem~\ref{th11}.}
If $\sX$ is not irreducible then it is decomposable and
we have $\sX = \sX^{(1)}\vee\sX^{(2)}$
where each $\sX^{(i)}$ satisfies $\sX^{(i)} = \psi (\sX^{(i)})$.
Lemma~\ref{le22} showed that $\sX$ and $\sX^{(i)}$ are generated by 
elements in the bounded region
$\| \bx \|_{\he{A}} \le R$. So the total multiplicity of elements
of each $\sX^{(i)}$ in the region $\| \bx \|_{\he{A}} \le R$
is strictly smaller than that of elements of $\sX$ in the region $\| \bx 
\|_{\he{A}} \le R$.
If one of $\sX^{(i)}$ is not irreducible then we can further partition
it, and decrease the total multiplicity in $\| \bx \|_{\he{A}} \le R$
again. Since $\sX$ is discrete, this process will end in finitely many steps,
yielding
$$
     \sX = \bigvee_{i=1}^N \sX^{(i)}.
$$
To see that the partition is unique, to each element $\sX_i$ 
of the partition 
is  associated a unique 
directed cycle in $\sG ( \sX )$ of multiplicity one. But each
directed cycle in  $\sG ( \sX )$, counted with multiplicity one,
must appear in some element of each partition.
\eproof

Another consequence  of Theorem~\ref{th21} is a  finiteness
condition on discrete irreducible multiset familities.

\begin{lemma}\label{le23}
Let $\psi(\cdot)$ be an inflation functional equation with
data $(\he{A}, \sD_{ij} )$. Then for
each $p \ge 1$ there are at most finitely many discrete multiset families
(possibly none)
$\sX$ satisfying $\psi(\sX) = \sX$ that are irreducible and 
have a periodic
cycle in $\sG ( \sX )$ of minimal period $p$.
\end{lemma}
\Proof 
Any periodic cycle $\{ (\bx_0, i_0), (\bx_1,i_1), \ldots, 
(\bx_{p-1}, i_{p-1}) \}$ 
constructed in Theorem~\ref{th21} has the form
$$
\bx_k = \he{A} \bx_{k-1} + \bd_k , \quad 1 \le k \le p
$$
with $\bd_k \in \sD_{i_{k}, i_{k-1}}$ and
$(\bx_p, i_p) = (\bx_0, i_0)$. It follows that
$$
\bx_0 = \he{A}^p \bx_0 + \sum_{k=1}^p \he{A}^{p-k} \bd_k .
$$
Thus
\beql{233}
    \bx_0 = - ( \he{A}^p - \he{I} )^{-1} \left(
\sum_{k=1}^p \he{A}^{p-k} \bd_k \right) ,
\eeq
where the matrix $\he{A}^p - \he{I}$ is invertible since
$\he{A}$ is expanding.
There are only a finite set of choices for the digits 
$\{ \bd_1, \ldots, \bd_p \}$, so the number of choices for $\bx_0$
is finite. 
\eproof

We continue the  study of discrete multiset families $\sX$ with 
$\psi ( \sX ) = \sX$ that are irreducible. We can give some
further information on the multiplicities that occur in such
$\sX$, under the extra hypothesis that the subdivision matrix
is primitive. 

\begin{theorem}[Multiplicity Dichotomy]\label{th22}
Suppose that the multiset family $\sX$ is discrete and
satisfies the inflation
functional equation $\psi ( \sX ) = \sX$ with primitive subdivision
matrix. 
If $\sX= (X_1, X_2, ..., X_n)$ is indecomposable, then exactly
one of the following cases holds:
\begin{itemize}
\item[\rm (i)]
Every element of each multiset $X_i$ has multiplicity at most $p$, 
the period of the
generating cycle in the directed graph $\sG ( \sX )$.
\item[\rm (ii)]
For each multiset $X_i$ in $\sX$,
 the multiplicities $m_{X_i}( \bx )$ of $\bx \in X_i$ are unbounded.
\end{itemize}
\end{theorem}
\paragraph{Remarks.}
(1) For the general case of imprimitive $\he{S}$, 
one can prove a weaker multiplicity  dichotomy; it 
applies only to those  $X_i$ such that $i$ lies in the 
strongly connected
component of the underlying graph of $\he{S}$ that
contains  the generating cycle.

(2) Both cases of the Multiplicity Dichotomy occur, as shown
by examples. Example~\ref{ex72b} gives an occurrence
of case (i) in which all $X_i$ have bounded multiplicity,
with some multiplicities exceeding one.
 
\paragraph{Proof.}
The multiplicity of a vertex $(\bx, i)$ in the graph
$\sG(\sX)$ is $m(\bx, i) := m_{X_i}(\bx)$. 
Fix a point $(\bx_0,i_0)$ in the generating cycle. 
We claim: For any point $(\bx,i) \ne (\bx_0, i_0)$ in 
$\sX$
the multiplicity $m(\bx, i)$ is equal to the number of 
cycle-free directed paths 
from $(\bx_0, i_0) $ to $(\bx,i)$. 
This is proved by induction
on the length $k$ of the longest cycle-free directed path 
from $\bx_0$ to $\bx$, which is bounded by (\ref{216}).
The base case $k =1$ is immediate.
We have 
$$
    m(\bx, i) = \sum_{(\bx', j) \in P(\bx, i)} m(\bx', j)
$$
where $P(\bx, i)$ is the set of preimages of $(\bx, i)$ in $\sG(\sX)$, 
Assuming $(\bx, i) \ne (\bx_0, i_0)$ all preimages have shorter
longest cycle-free path, and the induction step follows.

	Now assume that some vertex $(\bx_1, j)$ has multiplicity 
$m(\bx_1, j) = q > p$.
Without loss of generality we may assume that $(\bx_1,j)$ and 
$(\bx_0, i_0)$ have the same
color $j= i_0$, for if not by the primitivity hypothesis we
can always find a descendant $(\bx_1', i_0)$ of $(\bx_1, j)$, 
and clearly $m(\bx_1',i_0) \geq m(\bx_1, j)) > p$. 
By the pigeonhole 
principle there exist two cycle-free directed paths from $(\bx_0, i_0)$ to 
$(\bx_1, i_0)$, whose
lengths are $L$ and $L'$ respectively with $L \equiv L' \wmod{p}$. Say
$L' = L +sp$ for some $s\geq 0$. We can then create a second directed path of
length $L'$ by first going around the cycle $s$ times in the beginning
and then following the path of length $L$.
We show that these two different directed paths force
the multiplicity of vertices to be unbounded.

Since each directed edge can be labeled by elements in $\sD_{ij},~1 \leq 
i,j\leq m$, we
label the two directed paths by
$$
    {\mathcal P}_1 = (\bd_1, \dots, \bd_L), \shsp  
    {\mathcal P}_2 = (\bd_1', \dots, \bd_L').
$$
The fact that $\bx_0$ and $\bx_1$ have the same color $i$ implies that
$\bd_1\in \sD_{ji}$, $\bd_L\in \sD_{ik}$ and 
$\bd_1'\in \sD_{j'i}$, $\bd_L'\in \sD_{ik'}$ for some $j,\,k,\,j',\,k'$.
Evaluating the two paths yields
\beql{w2.1}
    \bx_1 = \A^L\bx_0 + \sum_{j=1}^L \A^{j-1}\bd_j
    = \A^L\bx_0 + \sum_{j=1}^L \A^{j-1}\bd_j'.
\eeq 
So $ \sum_{j=1}^L \A^{j-1}\bd_j = \sum_{j=1}^L \A^{j-1}\bd_j'$.
This means that following the two paths ${\mathcal P}_1$ and ${\mathcal P}_2$ 
from $(\bx_1, i_0)$ will lead to
the same vertex $(\bx_2, i_0)$, which also has color $i_0$. This process can
be continued to obtain vertices $(\bx_k, i_0)$, $k\geq 1$, 
all of which have the
same color $i_0$.
Note that there are at least $2^k$ distinct directed paths from 
$(\bx_0, i_0)$ to $(\bx_k, i_0)$ 
as we may concatenate ${\mathcal P}_1$ and ${\mathcal P}_2$ 
in
any combination. These $2^k$ directed paths remain distinct after
removing the cycles in the initial segment. Hence
$m(\bx_k, i_0) \geq 2^k$. 
So the multiplicity of $X_{i_0}$ is unbounded. Now the primitivity
of the subdivision matrix implies that any vertex $(\bx, i_0)$ has
a descendant in each $X_i$.
It follows that the multiplicity function is unbounded on every
$X_i$.
\eproof

As we have shown, all irreducible discrete multiset families $\sX$
are generated iteratively from a cycle as ``seed''. However, the multiset
family generated from a given cycle need not be discrete.  
The following lemma gives a criterion  for discreteness.

\begin{lemma}\label{le24}
Let $\psi(\cdot)$ be
 an inflation functional equation with data $( \he{A}, \sD_{ij} )$, 
and let $\sS^{(0)}$ be a finite multiset family with the inclusion property
$\sS^{(0)} \subseteq \sS^{(1)} = \psi ( \sS^{(0)})$.
Set $R = R( \he{A},  \sD )$ equal to the constant in (\ref{214}).
If for some $k \ge 1$,
\beql{228}
\sS^{(k)} \cap \{ \bx \in \RR^d : \| \bx \|_{\he{A}} \le R \} =
\sS^{(k+1)} \cap \{ \bx \in \RR^d : \| \bx \|_{\he{A}} \le R \}
\eeq
counting multiplicities, then the limit
$$\sX := \lim_{k \to \infty} \sS^{(k)}$$
is a discrete multiset family that satisfies
$\psi ( \sX ) = \sX$.
Conversely, if  $\sX$ is discrete then $($\ref{228}$)$ holds for all 
sufficiently large $k$.
\end{lemma}
\Proof
The iterative scheme $\sS^{(k+1)} = \psi ( \sS^{(k)} )$ is 
geometrically expanding outside the set
$\{ \bx \in \RR^d : \| \bx \|_{\he{A}} \le R \}$,
i.e. (\ref{216}) gives
$$
\| \he{A} \bx + \bd \| \ge \rho \| \bx \| - C
\ge \frac{\rho^2 +1}{\rho+1} \| \bx \| .
$$
The condition \eqn{228} says that $\sS^{(k)}$ stabilizes inside
$\{ \bx \in \RR^d: \| \bx \|_{\he{A}} \le R \}$.
By induction on $j \ge 0$ one obtains that
$\sS^{(k+j)}$ stabilizes inside the domain
$\left\{ \bx \in \RR^d : \| \bx \|_{\he{A}} \le 
\left( \frac{\rho^2 +1}{\rho+1} \right)^j R \right\}$.
The limit family $\sX$ is discrete inside each of these domains,
hence is discrete.
\eproof

Lemma~\ref{le24} leads to a recursively enumerable procedure  which 
recognizes all irreducible discrete sets $\sX$.
More precisely, consider the class of 
inflation functional equations $\psi(\cdot)$ whose data $(\he{A}, \sD_{ij})$ is
drawn from a computable subfield $\KK$ of $\CC$, which we take to be
 one in which addition
and multiplication are computable, and one can effectively test equality or
inequality of field elements. Then all finite cycles $Y$ for such
$\psi(\cdot)$ will have elements in $\KK$.  
One can give a procedure which
tests whether a given seed $\sS^{(0)}$ generates an irreducible
discrete $X$, using the criterion of Lemma~\ref{le24}.
It will eventually halt for all such sets,
and certify that they are discrete. (One does not
need to compute $R$ or $||\cdot||_{\he{A}}$ exactly to guarantee
halting.) 
This procedure also detects all  non-irreducible sets,
by either finding a second cycle or a cycle with multiplicity.
However this procedure is not recursive, for given an input cycle $Y$ that
generates an irreducible $\sX$ which is 
not discrete, the 
procedure can run forever without halting.

The following lemma gives  a simple case when all irreducible solutions to
the inflation functional equation are discrete (if any exist).

\begin{lemma}\label{le25}
Suppose that the inflation functional equation $\psi(\cdot)$
has data $( \he{A}, \sD_{ij} )$, in which $\he{A}$ is
an expanding integer matrix and all vectors in $\sD_{ij}$
have rational entries. Then all
irreducible multiset families
generated by a periodic cycle are discrete. 
\end{lemma}
\Proof
All cycles $Y$ of period $p$ are generated by solutions of
the form \eqn{233}, and the hypothesis guarantee that
such a solution $\bx$ is a rational vector. Starting from 
$(\bx, i)$ one generates the orbit $Y$ and a seed $\sS^{(0)}$
which consists of rational vectors. All elements of $\sS^{(k)}$
are rational vectors with denominators dividing  $dD$,
where $d$ is the denominator of $\bx$ and $D$ is the
greatest common denominator of all the rationals in $\sD_{ij}$.
Thus there are only a finite number of possible choices for
rational vectors in the ball of radius $R = R( \he{A},  \sD )$ in
\eqn{214}. Thus the  criterion of Lemma~\ref{le24} applies to
conclude that if 
$\sX = \lim_{k \to \infty} \sS^{(k)}$ exists, then it is
necessarily discrete. 
\eproof 

Lemma~\ref{le25} does not assert existence of
irreducible multiset families, only that they are discrete
if they do exist. Example~\ref{no-discrete} in \S8
shows there are cases satisfying the hypotheses where
no discrete solutions exist.
At the other extreme, there are 
inflation functional equations of this type 
having infinitely
many different irreducible discrete multiset families, see
Example~\ref{ex72a}.

We conclude this section by raising two related  (unsolved) computational
problems: 

\vspace{2mm}
\noindent{\em Computational Problem} (i).  
Given a seed $\sS_0$ which generates an irreducible
discrete multiset family $\sX$ satisfying a given inflation functional
equation $\psi(\sX) = \sX$ and an index $i$, determine whether  
the multiplicities of points in $X_i$ are bounded or not.

\vspace{2mm}
\noindent{\em Computational Problem} (ii).  Given the same data as above, 
determine the maximum multiplicity of points in $X_i$.

%
%
%
%
\setcounter{section}{3}
\setcounter{equation}{0}

\section{Perron Eigenvalue Condition}

In this section we consider 
inflation functional equations having
a primitive   subdivision matrix $\he{S}$
so that the Perron eigenvalue $\lambda(\he{S})$ is defined.
Our main object in this section is to show that a
necessary condition for an inflation functional
equation to have a solution that is a 
Delone set family is that the Perron eigenvalue 
$\he{S}$ satisfies the Perron eigenvalue condition
\beql{303}
\lambda (\he{S}) = | \rdet ( \he{A} ) |,
\eeq
as given in Theorem~\ref{th13}. To achieve this we
prove several preliminary results. Some of them apply
more generally to multiset families
 $\sX$ which are weakly uniformly discrete in the following sense.

\begin{defi}~\label{de110}
{\em
A multiset family $\sX$ is {\em weakly uniformly discrete} if the
union 
$\bigvee_{i =1}^n X_i$ is a weakly uniformly discrete multiset.
}
\end{defi}

A finite union of uniformly discrete multisets is 
weakly uniformly discrete if and only if each multiset itself
is  weakly uniformly discrete. 
So a multiset family $\sX$ is {\em weakly uniformly discrete} 
if and only if each $X_i$ in it is weakly uniformly discrete.

We first consider an irreducible multiset
family $\sX$ satisfying $\sX = \psi(\sX)$. 
Let $\sS^{(0)}$ be the generating cycle for $\sX$ and set
$\sS^{(k+1)} = \psi(\sS^{(k)})$. Write
$ \sS^{(k)}=(S^{(k)}_1, S^{(k)}_2, \dots, S^{(k)}_n)$
and let $H_i^{(k)}$ denote the total multiplicity of
the multiset $\sS^{(k)}_i$, which is given by 
$$
     H_i^{(k)} := \sum_{\bx\in \sS^{(k)}} m_{S^{(k)}_i}(\bx).
$$

\begin{lemma}  \label{le33}
   Let $\sX$ be an irreducible discrete multiset family satisfying 
$\psi(\sX) = \sX$ for the
data $(\he A, \sD_{ij})$, with subdivision matrix $\he{S}$.
Let $\bH^{(k)} = [H_1^{(k)}, \dots, H_n^{(k)}]^T$.
Then $\bH^{(k)} = \he S^k \bH^{(0)}$.
\end{lemma}
\Proof   We prove that $\bH^{(k+1)} = \he S \bH^{(k)}$. 
Observe that 
$$
  S^{(k+1)}_i = \bigvee_{j=1}^n (S^{(k)}_j +\sD_{ij}).
$$
Hence
$$
        H_i^{(k+1)} =  \sum_{j=1}^m |\sD_{ij}|H_j^{(k)},
$$
yielding $\bH^{(k+1)} = \he S \bH^{(k)}$.
Therefore $\bH^{(k)} = \he{S}^k \bH^{(0)}$.
\eproof

\begin{theorem}   \label{thm34}
 Let $\sX$ be an irreducible discrete multiset family satisfying 
$\psi(\sX) = \sX$ 
for the data $(\he A, \sD_{ij})$, with primitive subdivision matrix 
$\he{S}= [ |\sD_{ij}|]$.
Let $B_1(\bo)$ be the unit ball and set
\beql{w3.1}
    M_i^{(k)} := \sum_{\bx \in \us X_i \cap \he{A}^k(B_1(\bo))} m_{X_i}(\bx).
\eeq
Then there exist constants $C_1, C_2 >0$ such that, for each $k \ge 1$, 
\begin{equation}  \label{3.3}
     C_1 \lambda(\he S)^k \leq M_i^{(k)} \leq  C_2 \lambda(\he S)^k.
\eeq
\end{theorem}
\Proof  We first establish the lower bound. Observe that for each
element $\bx\in \bigvee_{i=1}^n S^{(k)}_i$ 
there exist digits $\bd_1, \dots, \bd_k$
in the collection of digit sets $\{\sD_{ij}\}$ such that
$\bx =\tau_{\bd_1}\circ \cdots \circ\tau_{\bd_k}(\bx_0)$, where $\bx_0$ is
in the generating cycle and $\tau_{\bd}(\bx):=\he A \bx +\bd$. So
$$
     \bx  = \he{A}^k\bx_0 + \sum_{j=1}^k \he A^{j-1} \bd_j
          = \he{A}^k\Bigl(\bx_0 + \sum_{j=1}^k ({\he{A}}^{-1})^{k+1 - j} \bd_j\Bigr).
$$
Since ${\he A}^{-1}$ is contracting, there exists a constant
$C$ depending only on $\{\sD_{ij}\}$, $\he A$ and $\sS^{(0)}$ such that
$$
    \Bigl\|\bx_0 + \sum_{j=1}^k ({\he A}^{-1})^{k+1-j} \bd_j\Bigr\| \leq C.
$$
Fix an integer $k_1>0$ so that $B_C(\bo) \subseteq \he{A}^{k_1}(B_1(\bo))$.  
Then $\bx \in \he{A}^{k+k_1}(B_1(\bo))$. It follows that for any 
$\bx\in \us{X}_i\cap \he{A}^{k+k_1}(B_1(\bo))$
we always have
$$
   m_{X_i}(\bx) = m_{S_i^{(k)}}(\bx).
$$  
Therefore $M_i^{(k+k_1)} \geq H^{(k)}$ for all $k \ge 1$.
By Lemma \ref{le33}
and the primitivity of $\he S$ we  have $H^{(k)} \geq C_1'\lambda(\he S)^k$
for some constant $C_1'$. The lower bound 
$M^{(k)} \ge C_1\lambda(\he S)^k$
follows by 
taking $C_1 = C_1'\lambda(\he S)^{-k_1}$.

To establish the upper bound we claim that there exists a
$k_2$ such that for each $k \ge 1$, each $i$ and 
$\bx\in \us X_i\cap \he A^{k}(B_1(\bo))$
and any vertex $\bx\in \us X_i$ with 
$\bx\in \he A^{k}(B_1(\bo))$ we must have
\beql{eq4.4}
        m_{X_i}(\bx) = m_{S_i^{(k+k_2)}}(\bx).
\eeq
The claim gives the upper bound, for  \eqn{eq4.4} 
implies  that $M_i^{(k)} \leq H^{(k+k_2)}_i$
holds for each $k \ge 1$ and each $i$.
However Lemma \ref{le33} and the primitivity of $\he S$ 
together show that
there exists a constant $C_2'$ such that 
$$M_i^{(k)} \leq C_2' \lambda(\he S)^{k+k_2}.$$
Combining these inequalities gives 
the upper bound $M^{(k)} \le C_2\lambda(\he S)^k$, 
on taking $C_2 = C_2' \lambda(\he S)^{k_2}$.

It remains to  prove the claim \eqn{eq4.4}. We note that any 
$\bx\in \us X_i\cap \he A^{k}(B_1)$ must be in $S_i^{(l)}$ 
for some $l$ as a result
of Lemma \ref{le21}. Suppose that $l \geq k$. 
Then there exist some digits $\bd_1, \bd_2, \dots, \bd_l$
in $\sD_{ij}$'s such that 
$$
     \bx  = \he{A}^k\bx_0 + \sum_{j=1}^l \he A^{j-1} \bd_j,
$$
and 
\beql{w3.2}
    \bx_0 = \he{A}^{-l}\bx - \sum_{j=1}^k \he A^{-j} \bd_{k+1-j}.
\eeq
Observe that $\bx_0$ is bounded since $l\geq k$, $\he A^{-k}\bx \in B_1$ and 
$\he{A}$ is expanding. Hence there exists a constant $R_0$ independent of
$k$ and $l$ such that
$\|\bx_0\| \leq R_0$. Let 
$$
   k_2 = \sum_{i=1}^n |\us X_i\cap B_{R_0}({\bf 0})|.
$$

We now argue by contradiction, 
and assume that (\ref{eq4.4}) is false. Then there exist some 
$k \ge 1$, $i$ and 
$\bx\in \us X_i\cap \he A^{k}(B_1)$ such that
$m_{X_i}(\bx) > m_{S_i^{(k+k_2)}}(\bx)$. This means there exists
a cycle-free directed path in $\sG(\sX)$ of length $l>k+k_2$
from an element $\bx_0$ in the generating cycle $\sS^{(0)}$ to $\bx$, say  
$\bx =\tau_{\bd_l}\circ \cdots \circ\tau_{\bd_2}\circ\tau_{\bd_1}(\bx_0)$.
Let $\bx_j=\tau_{\bd_{j}}(\bx_{j-1})$, $1\leq j \leq l$, be the vertices
of this directed path.
It follows that $\|\bx_j\| \leq R_0$ for $0 \leq j \leq l-k$.
But there are only $k_2$ vertices in $B_{R_0}({\bf 0})$ and $l-k >k_2$,
hence there exist two identical vertices among $\{\bx_j: ~1\leq j \leq l-k\}$.
This contradicts the assumption that the directed path is cycle-free, proving
(\ref{eq4.4}).
\eproof

\begin{theorem}   \label{thm35}
 Let $\sX$ be a weakly uniformly discrete $n$-multiset family satisfying 
$\psi(\sX) = \sX$ for the data $(\he A, \sD_{ij})$, 
with primitive subdivision matrix $\he S$.Then 
$$
\lambda (\he S) \leq |\det(\he A)|. 
$$
\end{theorem}
\Proof  We first assume that $\sX$ is irreducible.
The weak uniform discreteness of $\sX$ implies that 
$M_i^{(k)} \leq C {\rm Vol}(\he A^k B_1({\mathbf 0}))= C |\det(\he A)|^k$ 
for some positive constant $C$.
 From Theorem \ref{thm34} it
immediately follows that $|\det(\he A)| \geq \lambda(\he S)$.

If $\sX$ is reducible, then $\sX = \bigvee_{j=1}^N \sX^{(j)}$. The  
argument above now
applies to $\sX^{(1)}$ to yield $\lambda(\he S) \le |\det(\he A)|$.
\eproof

\medskip
\noindent
{\bf Proof of Theorem~\ref{th13}.}~~Since $\sX$ is a 
weak substitution Delone multiset
family, it is weakly uniformly discrete. By Theorem~\ref{thm35}
this  immediately gives
$\lambda(\he S) \le |\det(\he{A})|$. 

To prove the other direction
$\lambda(\he S) \ge |\det(\he{A})|$,
let $\sX = \bigvee_{j=1}^N \sX^{(j)}$.
Let $M_{i,j}^{(k)}$ be as in (\ref{w3.1}), but defined for $\sX^{(j)}$.
Then the relative denseness of each $X_i$  yields
$$
    \sum_{j=1}^N M_{i,j}^{(k)} \geq 
C' {\rm Vol}(\he{A}^k B_1({\mathbf 0})) = C' |\det(\he{A})|^k
$$
for some positive constant $C'$. Hence
$\max_j M_{i,j}^{(k)} \geq \frac{1}{N} C' |\det(\he{A})|^k$.
Taking the $k$-th roots and
letting $k \to \infty$ now yields
$\lambda(\he S)  \geq \ |\det(\he A)|$, using Theorem~\ref{thm34}.
\eproof

%
%
%
%
%

\setcounter{section}{4}
\setcounter{equation}{0}

\section{Multi-Tiling Functional Equation}

Our object in the next section is to relate 
the existence of weak substitution Delone multiset families to
self-replicating tilings. 
In this section we obtain preliminary information, 
concerning solutions of the multi-tiling functional
equation. Theorem~\ref{th51a} gives a necessary condition for
solutions of positive Lebesgue measure and 
Theorem~\ref{th44} 
gives a set of equivalent conditions for
 positive Lebsegue measure that apply
when the Perron eigenvalue condition  $\lambda({\he S})= |\det (\he{A})|$
holds.

We begin with  a basic existence result on
solutions to the multi-tiling equation. 

\begin{prop}\label{pr41}
The multi-tile functional equation
\beql{402}
\he{A} (T_i) = \bigcup_{j=1}^{n} (T_j + \sD_{ji} ) , \quad 1 \le i \le n ,
\eeq
with primitive subdivision matrix $\he{S}$ has a unique nonempty solution
$(T_1, \ldots, T_n )$ in which each $T_i$ is compact. 
In this solution  all $T_i$ are nonempty, and
\beql{w4.1}
  T_i = \Bigl\{\sum_{k=1}^\infty \A^{-k}\bd_{j_kj_{k-1}} \Bigm|
        \bd_{j_kj_{k-1}} \in \sD_{j_kj_{k-1}},
        ~(j_0,j_1,j_2,\dots) \in \{1,2,\dots,N\}^{\ZZ^+},~j_0 = i\Bigr\}.
\eeq
\end{prop}
\Proof 
Flaherty and Wang \cite[Proposition 2.3]{FW} 
prove under the hypothesis 
$$ 
(*)~~~\bigcup_{j=1}^n\sD_{ji} ~~~\mbox{ is nonempty for}~~~ 1 \le i \le n,
$$  
that the multi-tile equation (\ref{402}) has a unique nonempty solution 
$(T_1, \dots, T_n)$
in which {\em all} $T_i$ are compact sets, and that in this
solution all $T_i$
are nonempty. 
The primitivity assumption on $\he{S}$ implies that hypothesis $(*)$ holds, 
hence this result applies.
\eproof

We remark that (\ref{w4.1}) has a graph-theoretic
interpretation. Form a directed graph $\sG(\psi)$ whose vertices correspond
to sets $T_i$ and with a directed edge from $T_i$ to $T_j$ labelled
$\bd$ if $\bd \in \sD_{ji}$. Then \eqn{w4.1}  
is equivalent to saying that for 
each point $\bx\in T_i$
there exists an infinite directed path $(\bd_1,\bd_2,\bd_3,\dots)$ 
in the graph $\sG(\psi)$ where an edge
connects with $\bd_1 \in \sD_{ji}$ for some $j$ such that
\beql{w4.2}
   \bx = \sum_{k=1}^\infty \A^{-k}\bd_{k}
\eeq
and vice versa. 

We next give a simple necessary condition for
the sets $T_i$ to have positive Lebesgue measure.

\begin{theorem}~\label{th51a}
Suppose that the  multi-tile functional equation
\beql{400a}
\he{A} (T_i) = \bigcup_{j=1}^{n} (T_j + \sD_{ji} ) , \quad 1 \le i \le n ,
\eeq
with primitive subdivision matrix $\he{S}= [|\sD_{ij}|]$
has a nonempty solution $(T_1, ..., T_n)$ consisting of compact sets, in
which at least one $T_i$ has positive Lebesgue measure.
Then
\beql{400b}
\lambda(\he{S}) \ge |\det (\he{A})|.
\eeq
\end{theorem}

\Proof
Taking the Lebesgue measure of both sides of the multi-tile
equation \eqn{400a} gives, for $1 \le i \le n$, 
$$
    |\det (\he{A})| \mu(T_i) \leq \sum_{j=1}^n |\sD_{ji}| \mu(T_j),
$$
where $|\sD_{ji}|$ denotes the cardinality of $\sD_{ji}$ 
(counting multiplicity). By hypothesis  
${\mathbf v} := [\mu(T_1), \dots, \mu(T_n)]$ is a 
nonnegative row vector, not identially zero,
and the equation above gives
$${\mathbf v}\he{S} \geq |\det (\he{A})| {\mathbf v}.$$
Write ${\mathbf v}\he{S}= |\det (\he{A})| {\mathbf v} +{\mathbf w}_1  $,
where ${\mathbf w}_1$ is nonnegative. Repeated multiplication by $\he{S}$ and
back-substitution yields
${\mathbf v}\he{S}^k = |\det (\he{A})|^k {\mathbf v} + {\mathbf w}_k,$
where ${\mathbf w}_k $ is nonnegative. We conclude that the
spectral radius $\rho(\he{S})$ satisfies
$$
\rho(\he{S}) \ge 
\lim_{k \to \infty} ||\bv \he{S}^k||^{1/k} \ge |\det (\he{A})|.
$$
Thus $\lambda(\he{S}) \geq |\det (\he{A})|.$
\eproof

In the remainder of this section we develop criteria for positive
Lebesgue measures of the sets $T_i$ that apply when the
Perron eigenvalue condition $\lambda({\he S}) = |\det (\he{A})|$
holds.

We define
 the {\it digit multisets} 
\beql{403}
    \sD_{ji}^{m} :=
        \bigvee_{j_1, \ldots, j_{m-1} =1}^n
   ( \sD_{jj_{m-1}} + \he{A} \sD_{j_{m-1}  j_{m-2}} + \cdots + 
\he{A}^{m-1} \sD_{j_{1}i} ),
\eeq
in which the sum is interpreted as counting multiplicities.
 It is easy to check that iterating (\ref{402}) yields
\begin{equation} \label{404}
  {\he A}^{m}(T_{i})=\bigcup_{j=1}^{n} (T_{j}+\us{\cal D}_{ji}^{m}), ~\, 
i=1,\ldots,n.
\end{equation}
In \eqn{404} we do not count multiplicity, so it suffices to use
$\us {\cal D}_{ji}^{m}$ instead of ${\cal D}_{ji}^{m}$.

\begin{defi}
{\em 
    (i) A family of discrete multisets $\{{\mathcal E}_\alpha:~\alpha\in I\}$ 
in $\RR^d$
is {\em equi-uniformly discrete} if there exists an $\ep_0 >0$ such that
each ${\mathcal E}_\alpha$ is uniformly discrete and any two 
distinct elements in
${\mathcal E}_\alpha$ are separated by at least distance $\ep_0$. 
In particular, each ${\mathcal E}_\alpha$ is an ordinary set.

(ii) A family of discrete multisets 
$\{{\mathcal E}_\alpha:~\alpha\in I\}$ is called 
{\em weakly equi-uniformly discrete}
if there exists a fixed  $M>0$ such that for each $\alpha\in I$ 
and each ball $B$ of radius 1 in $\RR^d$ 
the number of elements of ${\mathcal E}_\alpha$ in $B$ 
(counting multiplicity) is bounded by $M$.
}  
\end{defi}

The following theorem is an extension of a theorem in 
Sirvent and Wang \cite{sirvent-wang}. Its  hypotheses can only
be satisfied when the Perron eigenvalue condition
holds, because the weak equi-uniformly discrete hypothesis
can only hold when $\lambda({\he S}) \le |\det (\he{A})|$

\begin{theorem} \label{thm42}
  Suppose that the
compact sets $(T_1,\dots, T_n)$ satisfy the multi-tile functional
equation
$$
    {\he A}(T_{i})=\bigcup_{j=1}^{n} (T_{j}+{\cal D}_{ji}), ~\, i=1,\ldots,n,
$$
with primitive subdivision matrix $\he{S}$, and that at least 
one $T_{i}$ has positive Lebesgue measure. 
If the collection of multisets 
$\{\us\sD_{ji}^{m}:~1\leq i,j \leq n, m \ge 1\}$ are 
weakly equi-uniformly discrete, then
each $T_{i}$ has nonempty interior  $T_i^{\circ}$, and is the
closure of its interior,  $T_{i} = \overline{T_i^{\circ}}$.
\end{theorem}
 
   The primitivity of the subdivision matrix implies
that if a single $T_i$ has positive Lebesgue measure, then they all do.
The proof of Theorem~\ref{thm42} is based on  the following 
covering lemma.

\begin{lemma} \label{le43}
 Suppose that $(T_1, ... , T_m)$ are compact sets
satisfying  an inflation
functional equation with primitive subdvision matrix,
and that at least one $T_i$ has
positive Lebesgue measure. 
Suppose also that the collection of multisets 
$\{\us\sD_{ji}^{m}:~1\leq i,j \leq n, m \ge 1\}$ is 
weakly equi-uniformly discrete.
Then given any sequence of positive numbers $\{\delta_{m}\}$
with $\lim_{m \to \infty} \delta_m = 0$,
there exist positive constants $R_{0}$ and $K_0$ such that 
the following holds:
For each $m \geq 1$
there exist subsets 
${\cal E}_1^{m}, \dots , {\cal E}_n^m$ of $\RR^d$ 
contained in the ball 
$B_{R_0}(\vec 0)$, 
each of cardinality bounded by $K_0$, such that
$\Omega_m:= \bigcup_{j=1}^n (T_j+{\cal E}_j^m)$ has the property that
\begin{equation} \label{405}
   \mu(B_{1}(\vec 0) \cap \Omega_m  )\geq 
(1-5^{d+1}\delta_{m}) \mu(B_{1}(\vec 0)).
\end{equation}
\end{lemma} 

\Proof 
Without loss of generality we assume $T_{1}$ has positive Lebesgue 
measure, so $T_1$ has a Lebesgue
point $\vec x^*$, i.e. there is a sequence $r_{m}\rightarrow 0$ such that
$$
 \mu( B_{r_{m}}(\vec x^*)  \cap    T_{1}    )\geq 
(1-\delta_{m})\mu(B_{r_{m}}(\vec x^*)).
$$
It follows that
\begin{equation} \label{406}
  \mu\Bigl( {\he A}^{l}(B_{r_{m}}(\vec x^*) \cap T_{1}    )\Bigr)\geq 
   (1-\delta_{m})\mu({\he A}^{l}(B_{r_{m}}(\vec x^*))),\, \mbox{ for all } 
l\geq 
0.
\end{equation}
We first show that for sufficiently large $l=l(m)$, there exists a unit ball
$B_{1}(\vec y)\subset {\he A}^{l}(B_{r_{m}}(\vec x^*))$ with
\begin{equation} \label{407}
       \mu\Bigl( B_{1}(\vec y) \cap  {\he A}^{l}(T_{1}      )\Bigr)\geq 
               (1-5^{d+1}\delta_{m})\mu(B_{1}(\vec 0)).
\end{equation}
Indeed, since ${\he A}$ is expanding,  
${\he A}^{l}(B_{r_{m}}(\vec x^*))$ is an 
ellipsoid
$O_{l,m}$ whose shortest axis goes to infinity as $l$ goes to infinity.
Let $O_{l,m}'$ be the homothetically shrunk ellipsoid with shortest axis 
decreased in length by $2$, so that all points in it are at distance at least 
1 from the boundary of $O_{l,m}$. By a standard covering lemma 
(cf. Stein~\cite[p. 9]{St70}) applied to
$O_{l,m}'$ there is a set  $\{B_{1}(\vec y')\}$ of disjoint unit balls
with centers in $O_{l,m}'$ that cover volume at least $5^{-d}\mu(O_{l,m}')$.
Also $5^{-d}\mu(O_{l,m}')\geq 5^{-d-1}\mu(O_{l,m})$ once the shortest axis
is of length at least $2(d+1)$. All these balls lie inside $O_{l,m}$.
By (\ref{406}) at most
$\delta_{m}\mu({\he A}^{l}(B_{r_{m}}(\vec x^{*})))$ of the volume of 
${\he A}^{l}(B_{r_{m}}(\vec x^*))$ is uncovered by 
${\he A}^{l}(B_{r_{m}}(\vec x^*)\cap T_{1})$,
so at least one of the disjoint balls $\{B_{1}(\vec y')\}$ must 
satisfy (\ref{407}).

   By (\ref{403}) we can rewrite the inequality
(\ref{407}) as 
$$
\mu\Bigl( B_{1}(\vec y) \cap 
\Bigl(\bigcup_{j=1}^n(T_{j}+\us{\cal D}_{j1}^{l})\Bigr) 
\Bigr)\geq 
  (1-5^{d+1}\delta_{m})\mu(B_{1}(\vec 0)),
$$
whence
$$
\mu\left( B_{1}(\vec 0) \cap 
\Bigl(\bigcup_{j=1}(T_{j}+\us{\cal D}_{j1}^{l}- \vec y)\Bigr)
\right) \geq 
  (1-5^{d+1}\delta_{m})\mu(B_{1}(\vec 0)).
$$
This shows that if we choose
$$
{\cal E}_j^{m} :=\Bigl\{\vec d- \vec y \, | \, \vec d\in\us{\cal D}_{j1}^{l} 
           \mbox{ with } (T_{j}+\vec d-\vec y)\cap B_{1}(\vec 0)\neq\emptyset 
\Bigr\}
$$
then (\ref{405}) holds. 
Since all $T_{j}$ are compact, all ${\cal E}_j^{m}$ lie  inside the ball 
$B_{R_{0}}(\vec 0)$ with $R_0:= 1+ \max_{i} {\rm diam}\,(T_i)$. 
The cardinality of
all ${\cal E}_j^{m}$ are upper bounded by  a constant $K_0$ because
the set of all  $\us\sD_{ij}^m$ is
weakly equi-uniformly discrete. 
\eproof

\vspace{3mm}

\noindent
{\bf Proof of Theorem~\ref{thm42}.}~~We apply Lemma~\ref{le43}
and choose a subsequence $m_k \to \infty$
 so that $\{{\cal E}_j^{m_k}\}$ converges
(as a multiset) for all $j$, and we denote the limit multisets 
by ${\cal E}_j^{\infty}$.
This can always be done because $\{{\cal E}_j^{m}\}$ are uniformly bounded
and have uniformly bounded cardinality. 
Clearly ${\cal E}_j^{\infty}$ has cardinality at most $K_{0}$.  So
\begin{eqnarray*}
   \mu\Bigl(B_{1}(\vec 0)\cap \Bigl(\bigcup_{j=1}^n
         (T_{j}+{\cal E}_j^{\infty})\Bigr)\Bigr) &\geq &
             \liminf_{k\rightarrow\infty} 
                 \mu\Bigl(B_{1}(\vec 0)\cap\Bigl(\bigcup_{j=1}^n
         (T_{j}+{\cal E}_j^{m_{k}})\Bigr)\Bigr)\\
  &\geq &\liminf_{k\rightarrow\infty} 
       (1-5^{d+1}\delta_{m_{k}})\mu(B_{1}(\vec 0))\\
  &=& \mu(B_{1}(\vec 0)).
\end{eqnarray*}
Since each $T_{j}+{\cal E}_j^{\infty}$ is a closed set, we must have
$$
     B_{1}(\vec 0)\cap\Bigl(\bigcup_{j=1}^n(T_{j}+{\cal E}_j^{\infty})
             \Bigr)=B_{1}(\vec 0).
$$
This means that at least one of the $T_j$'s 
must have nonempty interior $T_j^{\circ}$. But if so then
the primitivity of the subdivision matrix implies that all $T_j$ have 
nonempty
interior. Let $T_j' = \overline{T^{\circ}_j}$. Then $(T_1', \dots, T_n')$
must also satisfy the same multi-tile equation \eqn{402}. By the uniqueness
assertion in  Proposition~\ref{pr41} we have
$T_j = T_j'$ for $1 \le j \le n$. 
\eproof

We obtain the following 
characterization of the positive
measure property $\mu (T_i) > 0$, which applies when
the Perron eigenvalue condition holds.

\begin{theorem}\label{th44}
Assume that the family of compact sets $(T_1, \dots, T_n)$
satisfies a multi-tile functional equation with the data
$(\he A, \sD_{ij})$ with primitive subdivision
matrix $\he{S}$ satisfying the Perron eigenvalue condition
$\lambda(\he{S}) = |\det(\he{A})|.$
Then the following conditions are equivalent:
\begin{itemize}
\item[\rm (i)]
	For $m \ge 1$ all   multisets $\sD_{jk}^{m}$ for $1\leq i,j\leq n$ 
are ordinary sets, and the family 
$\{ \sD_{jk}^{m}:~1\leq i,j\leq n,\, m\geq1 \}$ 
is equi-uniformly discrete.

\item[\rm (ii)]
	For some fixed $1 \leq k,l \leq n$ and $m \ge 1$
the multisets $\sD_{kl}^{m}$ are ordinary sets and the family
$\{\sD_{kl}^{m}:~m\geq 1\}$
is equi-uniformly discrete.

\item[\rm (iii)]
	The family of multisets $\{\sD_{jk}^{m}:~1\leq i,j\leq n,\,m\geq1\}$ 
is  weakly equi-uniformly discrete.

\item[\rm (iv)]
	For some fixed $1 \leq k,l \leq n$ the family of multisets 
$\{\sD_{kl}^{m}:~m\geq 1\}$ is 
       weakly equi-uniformly discrete.

\item[\rm (v)]
	One set $T_j$ has $\mu (T_j) > 0$.

\item[\rm (vi)]
	Every set $T_j$ has $T_j = \overline{T_j^{\circ}}$, hence all 
$\mu(T_j) > 0$.

\item[\rm (vii)]
	Each $T_j = \overline{T_j^{\circ}}$, for $1 \le j \le n$,  and 
 $\mu ( \partial T_k) =0$.
\end{itemize}
\end{theorem}
\Proof
The implications (i) $\Rightarrow$ (ii)  $\Rightarrow$ (iv)
and  (i) $\Rightarrow$ (iii) $\Rightarrow$ (iv) are obvious. 

\medskip
\noindent
(iv) $\Rightarrow$ (v).~~First it is easy to check that $\he S^m = 
[|\sD_{ij}^m|]$
where $|\sD_{ij}^m|$ denotes the cardinality (counting multiplicity) of 
$\sD_{ij}^m$.
The primitivity hypothesis implies that there exists a $c_0>0$ such that,
for all $m \ge 1$,  
$|\sD_{ij}^{m}| \geq c_0 \lambda(\he S)^m$ 
for all $1\leq i, j \leq n$. 
Now for all $i$ let
$T^{(0)}_i = \overline{B_{1}(\bf 0)}$  and
$$
     T_i^{(m)} = \bigcup_{j=1}^n \he A^{-1}(T_j^{(m-1)} +\sD_{ji}), ~~m>0.
$$
Then $T_i^{(m)} \longrightarrow T_i$ in the Hausdorff metric (see \cite{FW}). 
We prove that $\mu(T_{l}) >0$ for all $1\leq l \leq n$. 
To see this, we note that
$$
     \he A^m \,(T_{i}^{(m)}) = \bigcup_{j=1}^n ~(T_j^{(0)} +\us\sD_{ji}^m\,).
$$ 
So
\begin{equation}  \label{408}
     \he A^m\,(T_{l}^{(m)}) = \bigcup_{j=1}^n~(T_j^{(0)} +\us\sD_{jl}^m)
	\supseteq T_{k}^{(0)} + \us\sD_{kl}^m.
\eeq
Since $\{\sD_{kl}^m\}$ are weakly equi-uniformly discrete, 
there exists a constant
$c_1>0$ such that 
$|\us\sD_{kl}^m| \geq c_1 |\sD_{kl}^m|$ for all $m$. 
Now since $\us\sD_{kl}^m$ are weakly equi-uniformly discrete there exists
an $M>0$ such that  each unit ball contains at
most $M$ elements of $\us\sD_{kl}^m$. Therefore each point in 
$ T_{k}^{(0)} + \us\sD_{kl}^m$ can be covered by no more than $M$ copies of
$T_k^{(0)}+\bd$, $\bd\in \us\sD_{kl}^m$.  Therefore by \eqn{408}
$$
   |\det(\he A)|^m \mu(T_{l}^{(m)}) 
	\geq \frac{1}{M}|\us\sD_{kl}^m| \mu(T_{l}^{(0)}) 
	\geq \frac{1}{M}\,\delta\,c_0 \,c_1\, \lambda(\he S)^m,
$$
where $\delta$ is the volume of the ball $B_{1}(\bf 0)$. 
Using the hypothesis $\lambda(\he{S})= |\det(\he A)|$ yields 
$$
\mu(T_{l}^{(m)}) \geq \frac{1}{M}\,\delta\,c_0\,c_1>0.
$$ 
It follows that, for all $l$,
$$
    \mu(T_{l}) \geq \liminf_{m\ra\infty} \mu(T_{l}^{(m)}) >0.
$$

\medskip
\noindent
(v) $\Rightarrow$ (i).~~First we note that $\mu(T_j) >0$ for all $j$ as a 
result of
the primitivity of $\he S$. Let $\vec e = [\mu(T_1), \mu(T_2), \dots, 
\mu(T_n)]$.
Taking the Lebesgue measure of both sides of the 
iterated multi-tile equation \eqn{404} yields
\begin{equation}  \label{409}
     \lambda^m(\he S) \mu(T_i) \leq \sum_{j=1}^n \mu(T_j)|\us\sD_{ji}^m|  
     \leq \sum_{j=1}^n \mu(T_j)|\sD_{ji}^m|,~~1\leq i \leq n.
\eeq
In other words, $\lambda^m(\he S) \vec e \leq \vec e\he S^m$. 
But $\he S^m$ is a primitive
nonnegative matrix with Perron-Frobenius eigenvalue $\lambda^m(\he S)$. 
So \eqn{409} can hold only when $\vec e$ is a 
left Perron-Frobenius eigenvector of $\he S$ 
and all inequalities in \eqn{409} are equalities. This immediately 
yields
$|\sD_{ji}^m| =|\us\sD_{ji}^m|$, hence each  $\sD_{ji}^m$ is an ordinary set.
The equalities in \eqn{409} also imply that the unions $T_j + \us\sD_{ji}^m$
are all measure-wise disjoint. So all $\us\sD_{ji}^m$, and hence all
$\sD_{ji}^m$, are equi-uniformly discrete for some $\ep>0$.

\medskip
\noindent
(v) $\Leftrightarrow$ (vi).  ~~This follows from Theorem~\ref{thm42},
 since (v) also implies (i).

\medskip
\noindent
(vi) $\Rightarrow$ (vii). 
~~We only need to prove that $\mu(\partial T_j)=0$ for 
all $j$. Let
$\vec v = [\mu(\partial T_1), \dots, \mu(\partial T_n)]$. We have for all $i$
$$
    \he A^m(\partial T_i) = \partial(\he A^m(T_i)) =
   \partial \Bigl(\bigcup_{j=1}^n (T_j +\us\sD_{ji}^m)\Bigr)
   \subseteq \bigcup_{j=1}^n (\partial T_j + \us\sD_{ji}^m).
$$
Similar to \eqn{409}, taking the Lebesgue measure yields
$\lambda^m(\he S) \vec v \leq \vec v \he S^m$. Again, this can occur only when
$\vec v = \bf 0$ or $\vec v$ is a Perron-Frobenius left eigenvector of $\he 
S^m$. 
Assume that $\vec v \neq \bf 0$ then all $\mu(\partial T_j) >0$, and
\begin{equation}  \label{410}
    \mu\Bigl(\partial(\he A^m(T_i))\Bigr) = 
    \sum_{j=1}^n |\sD_{ji}^m| \mu(\partial T_j).
\eeq
But $T_i$ has nonempty interior, so for sufficiently large $m>0$ the inflated
set $\he A^m(T_i)$ will contain a sufficiently large ball in its interior.
Since $\he A^m(T_i)$ is the union of $ T_j + \us\sD_{ji}^m$, $1 \leq j \leq n$,
there must be some $k$ and $ \bd\in\us\sD_{ki}$ such that $T_k+\bd$ is 
completely
contained in the interior of $\he A^m(T_i)$. Hence
$$
   \partial(\he A^m(T_i)) \subseteq  \bigcup_{j=1}^n (\partial T_j + 
\us\sD_{ji}^m)
	\setminus (T_k+\bd).
$$
So \eqn{410} is impossible, a contradiction.

\medskip
\noindent
(vii) $\Rightarrow$ (vi). ~~ This is obvious.
\eproof

%
%
%
\setcounter{section}{5}
\setcounter{equation}{0}

\section{ Weakly Uniformly Discrete Multiset Families and Tilings}

In this section we suppose that the Perron eigenvalue condition
$\lambda(\he{S}) = |\rdet(\he{A})|$
holds. Using the results of \S5 we relate the existence
of positive Lebesgue measure tiles $T_i$ for the associated
multi-tiling equation to weakly uniformly discrete
multiset solutions and self-replicating multi-tiling solutions 
to some iterate  $\psi^N(\cdot)$ of the inflation functional equation.
Then we  deduce Theorem~\ref{th14}. 

\begin{theorem}\label{th41}
Let $\psi ( \sX ) = \sX$ be an inflation functional equation 
that has a primitive subdivision matrix $\he{S}$ that
satisfies the Perron eigenvalue condition
$$\lambda(\he{S}) = |\rdet(\he{A})|.$$
Then the following conditions are equivalent:
\begin{itemize}
\item[\rm (i)]
For some $N > 0$ there exists a weakly uniformly discrete multiset family 
$\hat\sX$ such that $\psi^N (\hat \sX ) = \hat\sX$.
\item[\rm (ii)]
For some $N > 0$ there exists a uniformly discrete multiset family 
$\hat\sX$ such that $\psi^N (\hat \sX ) = \hat\sX$.
\item[\rm (iii)]
For some $N > 0$ there exists a 
self-replicating multi-tiling family $\hat{\sX}$ with
 $\psi^N(\hat \sX ) = \hat\sX$.
\item[\rm (iv)]
The unique compact solution $(T_1, \ldots, T_n )$ 
of the associated multi-tile functional equation 
consists of sets $T_i$ that have positive Lebesgue measure, 
$1 \le i \le n$.
\end{itemize}
\end{theorem}

\noindent
{\bf Proof.}~~Iterating $N$ times the 
inflation  functional
equation $\sX = \psi(\sX) $ on multisets
gives a new inflation functional equation
$\sX = \psi^N(\sX)$, which corresponds to
\begin{equation}  \label{411}
     X_i = \bigvee_{j=1}^n (\he A^N (X_j) + \sD_{ij}^{N}), 
\quad 1\leq i \leq n,
\end{equation}
where the sums are interpreted as multiset sums.

We show that  (i) $\Rightarrow$ (iv) $\Rightarrow$ (iii) $\Rightarrow$ (ii)
$\Rightarrow$ (i). To begin, the  implications 
(iii) $\Rightarrow$ (ii) $\Rightarrow$ (i) are obvious. 

\medskip
\noindent 
(i) $\Rightarrow$ (iv).~~By assumption there exists an $M>0$ such that
any unit ball in $\RR^d$ contains at most $M$ elements (counting multiplicity)
of each $X_i$.  Replacing $N$ by $mN$ in \eqn{411} 
it follows that any unit ball 
contains at
most $M$ elements of each $\sD_{ij}^{mN}$.  Observe that
by \eqn{403} we have
$$
     \sD_{ij}^{k} =  \bigvee_{l=1}^n (\sD_{il}^{k-1} + 
\he{A}^{k-1}\sD_{lj}).
$$
Therefore any unit ball contains at most $M$ elements of each 
$\sD_{ij}^{k-1}$ if
it contains at most $M$ elements of each $\sD_{ij}^{k}$.
This immediately yields the  weakly equi-uniform discreteness of all the sets  
$\{\sD_{ij}^{k}\}$.
So $T_i^o \neq \emptyset$ by Theorem~\ref{th44}.

\medskip
\noindent 
(iv) $\Rightarrow$ (iii).~~Since $T_1^o \neq \emptyset$ it follows from the
observation (\ref{w4.2}) that there exists an infinite directed path
$(\bd_1,\bd_2,\bd_3,\dots)$ in the graph
$\sG(\sX)$ with $\bd_1 \in \sD_{j1}$ for some $j$ such that
$$
   \bx_0 = \sum_{k=1}^\infty \A^{-k}\bd_{k} \in T_1^o.
$$
Since all $\sD_{ij}$ are bounded, there exists an $N'>0$ such that
for all infinite directed paths $(\bd_1',\bd_2',\bd_3',\dots)$
with $\bd_j' = \bd_j$ for $j \leq N'$ we also have
$$
   \bx_0' = \sum_{k=1}^\infty \A^{-k}\bd_{k}' \in T_1^o.
$$
The primitivity of the subdivision matrix $\he{S}$ now implies that 
we can find an infinite directed paths $(\bd_1^*,\bd_2^*,\bd_3^*,\dots)$
which has $\bd_j^* = \bd_j$ for $j \leq N'$, and which is periodic
for some period $N\geq N'$ in the sense that $\bd_{k+N}^* = \bd^*_{k}$
for all $k$. Let $\bx_0^* = \sum_{k=1}^\infty \A^{-k}\bd_{k}^*$. Then
$\bx_0^* \in T_1^o$ and $\A^N \bx_0^* = \bx_0^*+ \bd$ for
$\bd = \sum_{j=0}^{N-1} \A^j \bd_{N-j}^*$. 
Observe that $\bd\in \sD_{11}^N$. So we have
\beql{w4.3}
   -\bx_0^* \in \A^N (-\bx_0^*) + \sD_{11}^N.
\eeq

Now consider the inflation functional equation $\sX = \psi^N(\sX)$. 
Set $\sX^{(0)} = (X_1^{(0)}, \dots, X_n^{(0)})$, with
$$
    X_1^{(0)} = \{-\bx_0^*\},~X_2^{(0)} = \emptyset,~\dots,~ X_n^{(0)} = 
\emptyset.
$$
Define $\sX^{(m)} = \psi^N(\sX^{(m-1)}) =\psi^{mN}(\sX^{(0)})$.
Set
$$
   X_1^{(1)} = \bigvee_{j=1}^n (\he A^N (X_j^{(0)}) + \sD_{1j}^N) = 
\he{A}^N(\{\bx_0^*\}) + \sD_{11}^N,
$$
hence we have $X_1^{(0)} \subseteq  X_1^{(1)}$
by (\ref{w4.3}), and obviously we have 
$ \emptyset=  X_i^{(0)}\subseteq  X_i^{(1)}$ for $i\geq 2$.
So the inclusion property
$\sX^{(0)} \subseteq \psi^N(\sX^{(0)}) = \sX^{(1)}$ holds.
It follows that $\sX^{(0)} \subseteq \sX^{(1)}\subseteq \sX^{(2)} \subseteq 
\cdots$.
But notice that
$$
 X_i^{(m)} = \bigvee_{j=1}^n (\he A^{mN} (X_j^{(0)}) + \sD_{ij}^{mN}) 
	= \he A^{mN}\{\bx_0^*\} + \sD_{i1}^{mN}.
$$
We conclude that
 each $X_i^{(m)}$ is an ordinary set and is $\ep_0$-uniformly discrete. Let
$X_i = \bigcup_{m=0}^\infty X_i^{(m)}$ and $\hat\sX= (X_1, \dots, X_n)$.
Then $\hat\sX = \psi^N(\hat\sX)$, and $\hat\sX$ is $\ep_0$-uniformly
discrete.

It remains to show that $\hat\sX$ is  a Delone family and
to establish the tiling property of $\hat\sX$. Observe that
$\bf 0$ is in the interior of $T_1 -\bx_0^*$. Now
$$
  {\he A}^{mN}(T_{1}- \bx_0^*)=
    \bigcup_{j=1}^{n} (T_{j}+\us{\cal D}_{j1}^{mN}-\A^{mN}\bx_0^*)
	= \bigcup_{j=1}^{n} (T_{j} + X^{(m)}_j).
$$
Note that the union $\bigcup_{j=1}^{n} (T_{j}+\us{\cal D}_{j1}^{mN})$ is 
measure-wise
disjoint as proved in the proof of (v) $\Rightarrow$ (i) in
Theorem \ref{th44}.  
Taking the
limit as $m \rightarrow \infty$ will keep the union 
measure-wise disjoint.
Hence $\bigcup_{j=1}^{n} (T_{j} + X_j)$
is a tiling of $\RR^d$. So all $X_j$ must be relatively dense as a 
result of the
primitivity of $\he{S}$. This completes the proof.
\eproof

We conclude this section by  deducing Theorem~\ref{th14} as a consequence of
Theorem~\ref{th41}.

\noindent\paragraph{Proof of Theorem~\ref{th14}.}
\noindent
(iii) $\Rightarrow$ (ii). This follows from Theorem~\ref{th41} (iv)
$\Rightarrow$ (iii). 

\noindent
(ii) $\Rightarrow$ (i). A self-replicating multi-tiling family is a
weak substitution Delone multiset family.

\noindent
(i)$\Rightarrow$ (iii). By Theorem~\ref{th13} the inflation functional
equation satisfies the Perron eigenvalue condition. A 
weak substitution Delone multiset family 
is a weakly uniformly discrete multiset family, and satisfies the
Perron eigenvalue condition by Theorem~\ref{th13}.
Thus the conditions of 
Theorem~\ref{th41} (i) hold, and the result then
follows from Theorem~\ref{th41} (iv).
\eproof

%
%
%

\setcounter{section}{6}
\setcounter{equation}{0}

\section{Self-Replicating Multi-Tiling Families}

In this section we study self-replicating multi-tiling families as a 
subclass of weak substitution Delone multiset families.

\begin{theorem}\label{th51}
Let $\sX$ be an irreducible weak Delone multiset
family satisfying the inflation functional equation
$\psi(\sX) = \sX$ 
for the data $( \he{A}, \sD_{ij} )$, where the subdivision matrix $\he{S}$
is primitive.
Suppose that the fundamental cycle of $\sX$ has period $1$.
Then ${\sX}$ is a self-replicating multi-tiling family.
\end{theorem}

\Proof Since the fundamental cycle of $\sX=(X_1,..., X_n)$ 
has period 1,
it contains a single element $\{\bx_0\}$, and since it is
irreducible, by Theorem~\ref{th12}
it has multiplicity one. 
Without loss of generality we assume that
$\bx_0 \in X_1$, so $\bx_0 = \he A\bx_0 +\bd$ for some $\bd \in \sD_{11}$.
Let $\sS^{(0)} = (S_1^{(0)}, \dots, S_n^{(0)})$ such that
$S_1^{(0)} = \{\bx_0\}$ and all other $S_i^{(0)} = \emptyset$.
Define $\sS^{(m)} := \psi(\sS^{(0)}) = (S_1^{(m)}, \dots, S_n^{(m)})$.
It follows from the expression for $\psi^m$ given in (\ref{411}) that
\begin{equation} \label{5.1}
  S^{(m)}_i = \bigvee_{j=1}^n (\he A^m (S_j^{(0)}) + \sD_{ij}^{m})
	= \he A^m\bx_0 + \sD_{i1}^m.
\end{equation}
Suppose that $(T_1, \dots, T_n)$ is the set of self-affine multi-tiles
corresponding to $(\he A, \sD_{ij})$. Now by
hypothesis Theorem \ref{th14} (i)
 holds, so by Theorem \ref{th14} (iv) 
each $T_i$ has positive
Lebesgue measure, and
each $T_i$ satisfies $\overline {T_i^{\circ}} = T_i$. We have
\begin{equation} \label{5.2}
  \bigcup_{j=1}^n (T_j+ S^{(m)}_j) = 
    \bigcup_{j=1}^n (T_j+\he A^m \bx_0  + \sD_{j1}^{m}).
\end{equation}
It follows from 
$ \he A^m (T_1) = \bigcup_{j=1}^n (T_j+ \sD_{j1}^{m}) $ that 
\begin{equation} \label{5.3}
 \he A^m (T_1+\bx_0) = \bigcup_{j=1}^n (T_j+ S^{(m)}_j).
\end{equation}
The unions on the right side of (\ref{5.3}) are all measure-wise disjoint.
Taking the limit $m \ra \infty$ we see that
$\Omega := \bigcup_{j=1}^n (T_j +X_j)$ is a packing of $\RR^d$.

It remains to  prove that $\Omega$ is a tiling.
The set 
$\Omega = \bigcup_{j=1}^n (T_j +X_j)$ is closed and satisfies
\begin{eqnarray*}
   \he A (\Omega) &=& \bigcup_{j=1}^n \Bigl(\he A(T_j) +\he A(X_j)\Bigr) \\
   &=& \bigcup_{j=1}^n \Bigl(\bigcup_{i=1}^n (T_i+\sD_{ij} +
\he A(X_j))\Bigr) \\
   &=& \bigcup_{i=1}^n \Bigl(T_i + \bigcup_{j=1}^n 
(\sD_{ij} +\he A(X_j))\Bigr) \\
   &=& \bigcup_{i=1}^n (T_i + X_i) ~=~ \Omega. 
\end{eqnarray*}
We now  argue by contradiction, and suppose $\Omega \neq \RR^d.$
Since $\Omega$ is closed, there exists a ball $B$ of radius $\ep>0$ such that
$ B \cap \Omega = \emptyset$, which yields $\he A^N (B) \cap \he A^N(\Omega)
= \he{A}^N (B) \cap \Omega = \emptyset$. But $\he A$ is expanding.
So by taking $N$ sufficiently
large $\he A^N(B)$ contains a ball of arbitrarily large radius. This ball
is disjoint from $\Omega$, so it is not
filled by any translate $T_j + \bx_j$, $\bx_j \in X_j$ and $1 \leq j \leq n$.
Therefore $X_j$ cannot be a weak Delone set, a contradiction. Thus
we have a tiling.
\eproof

\paragraph{Remarks.} (1) There exist irreducible self-replicating 
multi-tiling families
whose fundamental cycles have period exceeeding $1$, see
Example \ref{ex-selfrep}. Thus the condition of Theorem~\ref{th51}
is only a sufficient condition (not a necessary
condition) to give a  self-replicating multi-tiling family.

\noindent
(2) Let  $\sX$ be an irreducible Delone set family which 
satisfies an inflation functional equation 
$\psi(\sX) = \sX$ with primitive subdivision matrix $\he{S}$
has a fundamental cycle of
period $p$. By Theorem~\ref{th13} and Theorem~\ref{th41}
the unique solution $\sT :=(T_1,..., T_n)$ consisting of compact
sets of the
associated multi-tile equation
has $T_i$ of positive measure, which are the closure
of their interiors.  
The set $\sT + \sX := \bigcup_{i = 1}^n (T_i + X_i)$ 
gives a partial $q$-packing of $\RR^d$ for some $q\leq p$, using the 
tiles $T_i$. That is, each point of $\RR^d$ is covered with
multiplicity at most $q$ off a set of measure zero, a set of infinite
measure has multiplicity exactly $q$, and possibly another set of infinite
measure has strictly smaller multiplicity. 

One expects that in many cases this construction produces a
$p$-thick  multi-tiling of $\RR^d$.
We formulate the following problem:
\vspace{3mm}

\noindent
{\it Problem.}  Let  $\sX$ be an irreducible Delone set family which 
satisfies an inflation functional equation 
$\psi(\sX) = \sX$ with primitive subdivision matrix $\he{S}$
and is uniformly discrete and has a
fundamental cycle of period $p$. Is it true that $\sT + \sX := \bigcup_{i = 1}^n 
(T_i + X_i)$
is always a $q$-thick multi-tiling for some $q\leq p$? If not, 
what are extra conditions needed to ensure it? 

%
%
%
%
\setcounter{section}{7}
\setcounter{equation}{0}

\section{Examples}
\setcounter{equation}{0}

\begin{exam}\label{ex71}(Substitution multiset 
with unbounded multiplicity function)
{\rm
Let $\he{A} = [3]$ and $\sD_{1,1} = \{ 0,1,2, 3 \}$, so
that $|\sD_{1,1}| = 4 > 3 = |\det(\he{A})|.$
Take the seed $\sS^{(0)} =  S_1^{(0)} = \{0\}$.
Then the inclusion property $\sS^{(0)} \subseteq \sS^{(1)}$ holds, hence
$$X = X_1 = \lim_{k \to \infty}  S_1^{(k)}$$
defines a multiset $X$.
The multiset $X \subseteq \ZZ_{\ge 0}$ and each  point 
$l \in \ZZ^+$ occurs with 
finite multiplicity $m_{X}(l)$. One can verify that 
$m_{X}(3^k) = k+1$, which shows that 
$m(l)$ is unbounded as $l \to \infty$.
This example corresponds to case (ii) in Theorem~\ref{th22}.
}
\end{exam}

\begin{exam}\label{ex72}(Discrete substitution set family
that is not uniformly discrete)
{\rm Let $\he{A} = [3]$ and $\sD_{1,1} = \{ 0,1, \pi \}$, with 
$\pi = 3.14159 \ldots$.
Take the seed $\sS^{(0)} =  S_1^{(0)} = \{0\}$.
Then $\sS^{(0)} \subseteq \sS^{(1)}$ and the limit
$$X_1 = \lim_{k \to \infty}  S_1^{(k)}$$
exists.
In this case the multiset $X_1 \subseteq \RR_{\ge 0}$ is discrete,
and its elements all have multiplicity one.
It is easy to show that it has linear growth.
Indeed the $2 \cdot 3^{n-1}$ elements in $ S_1^{(n)} \diagdown
 S_1^{(n-1)}$ all satisfy
$$3^n \le x \le \pi (3^n + 3^{n-1} + \cdots + 3+1) \le 2 \pi \cdot 3^n ~.$$
The associated multitile functional equation is
$$\he{A} (T_1) = T_1 \cup (T_1 + 1) \cup (T_1 + \pi ) ~.$$
The compact solution $T_1$ to this equation has Lebesgue measure zero,
see Kenyon \cite{Ke97} or Lagarias and Wang \cite{LW96d}.
It follows from Theorem~\ref{th41} that $X_1$ cannot be uniformly discrete.
}
\end{exam}

\begin{exam}~\label{ex72a}(Inflation functional equation having
infinitely many discrete solutions)
{\rm
Consider the inflation functional equation with $\he{A} = [2]$ on $\RR$
with $n=1$ and $\sD_{1,1} = \{0,1\}$,
This satisfies the hypotheses of Lemma~\ref{le25}. 
The associated tile is $T_1 = [0, 1].$ The allowed starting
points for a cycle $Y$ of
period $p$ (an arbitrary positive integer) are given by \eqn{233}, which gives
$$x_0 = -\frac{m}{2^p -1} \qquad\mbox{for}\qquad 0 \le m \le 2^p - 1.$$
Each such $x_0$ generates
an irreducible multiset $\sX_{m, p}$ consisting of a single
multiset $X_1 = X_1(m,p)$ which has all multiplicites  one.
Now suppose that $p \ge 2$ and restrict to
those  $x_0 = - m/(2^p - 1)$ is in a cycle of minimal period $p$,
with $p \ge 2$.
This always happens when $g.c.d.(m, 2^p -1) = 1$, and possibly in
other cases as well, but rules out the cases $m=0$ and $m = 2^p -1$, which
generate minimal cycles of period $1$, having the  symbolic dynamics
$d_1 d_2 \cdots d_p = 0^p$ and $1^p$, respectively.
Thus  $-1 < x_0 < 0$, 
so that the tile $T_1 + x_0$ includes $0$ in
its interior, and the same holds for the other $p-1$ tiles in
the periodic cycle.
We conclude that the cycles of minimal period $p$
generate  a multiple tiling of $\RR$ with thickness $p$, using
copies of the tile $T_1$.
}
\end{exam}

\begin{exam} \label{no-discrete}
(Inflation functional equation having no  
nonempty discrete multiset solution)
{\rm 
Consider the inflation functional equation with  $\he A = [2]$ in $\RR$,
and $n=2$ with digit sets 
$\sD_{11} = \sD_{12}=\sD_{21}=\sD_{22}=\{0,1\}$.
This data
satisfies the hypotheses of Lemma~\ref{le25}, but the
inflation functional equation $\psi(\sX) = \sX$ has no discrete 
multiset solution. If fact, all elements in
$\sX$ must have infinite multiplicity. To see this we iterate the 
inflation functional  equation to obtain
obtain 
$$
     X_i = (X_1 + \sD_{i1}^m) \vee (X_2+\sD_{i2}^m),~~i=1,2.
$$
Now each $\sD_{ij}^m$ has an underlying set 
$\us\sD_{ij}^m = \{0,1,\dots, 2^{m}-1\}$
with each element having multiplicity $2^m$. 
Therefore the multiplicity of each element in $X_i$ is at least $2^m$. 
Hence no element in any $X_i$ can have a finite multiplicity.
}
\end{exam}

\begin{exam}\label{ex72b}
(Irreducible weak substitution Delone multiset family with bounded
multiplicities, some multiplicities exceeding one)
{\rm
Consider the inflation functional equation
$\he{A} = [3]$ on $\RR$, with $n=2$ and
digit sets  $\sD_{1,1} = \{\pi + 3\}$,
$\sD_{1,2} = \{ 1\}, $ $\sD_{2, 1}= \{ -3, \pi \}$,
$\sD_{2,2} = \emptyset.$ The subdivision matrix 
$\he{S} = \left[ \begin{array}{rr}
            1 & 1 \\
            2 & 0 \end{array} \right]
$  
is primitive,
and its Perron eigenvalue $\lambda(\he{S}) = 2.$
We claim that the  cycle $Y = \{ 0 \in X_1, 1 \in X_2 \}$ 
of period $p = 2$, both points of multiplicity 1,  generates a
discrete irreducible multiset
family $\sX = (X_1, X_2)$.
The only point reached in
one step by exiting from the cycle is $ y =\pi + 3 \in X_1$,
which is reached in two different ways, one from $0 \in X_1$ and
one from $1 \in X_2$, so  
has multiplicity $2$. All other points are descendents
of $y$, and necessarily fall in the interval $[ \pi + 3, \infty)$.
All the maps $x \to 3x + d$ for $d \in \sD_{i,j}$
are expanding outside the interval
$[-4,4]$ with expansion factor at least $1.4$, so it follows that 
$\sX = (X_1, X_2)$ exists, and is irreducible and discrete.
We claim that the  descendents
of $y$ are distinct and consequently all have
multiplicity two. To show they are distinct,
note that all descendants of $y$  have the form 
$x= (3^k a_0 + 3^{k-1} a_{1} + ... +3 a_{k-1} + a_k)\pi + r,$
in which $a_0=1$ and each subsequent  $a_i = 0$ or $1$,
 and $r$ is an integer. 
The sequence $(a_0, a_1..., a_k)$ completely specifies the digit
sequence leading to $x$, which determines which
of $X_1$ and $X_2$ the point
$x$ belongs to; furthermore every such digit sequence is legal.
Since all integers  $3^k a_0 + 3^{k-1} a_{1} + ... +3 a_{k-1} + a_k$
are distinct, and $\pi$ is irrational, we conclude
that all $x$ are distinct. Thus all points in $X_1$ and $X_2$ have
multiplicity two, except the generating cycle $Y$, whose two points
have multiplicity one.
}
\end{exam}
\begin{exam}\label{ex73}(Substitution Delone set  
that is not self-replicating)
{\rm Let $\he{A} = [3]$ and $\sD_{1,1} = \{-1,0,1\}$.
The associated tile $T_1 = [ - \frac{1}{2} , \frac{1}{2} ]$, 
which consists of  all balanced  ternary expansions
$$x = \sum_{j=1}^\infty \bd_j 3^{-j}, \quad
\bd_j \in \{-1,0,1\} ~.
$$
The set
$\sS^{(0)} = \{ x_0 = - 1/8 , ~~x_1 = - 3/8 \}$
has 
$$\sS^{(0)} \subseteq \psi ( \sS^{(0)} ) =
\{ - 17/8 , -11/8 , -9/8 , -3/8 , -1/8 , 5/8 \},$$ 
hence
 $\sS^{(0)}$ generates an irreducible discrete multiset family
$\sX = \{ X_1 \}$ which consists of the single set $X_1$ given by
$$X_1 = \sum_{k \to \infty} \sS^{(k)} ~.$$
The set $X_1$ is irreducible and $\sS_0$ is its generating cycle of 
period $2$.
A calculation gives
$$X_1 = \left( - \frac{3}{8} + \ZZ \right) \cup \left( - \frac{1}{8} \cup \ZZ 
\right) ~.$$
It is a Delone set, and $X_1 + T$ is a {\em multiple tiling} of $\RR$ of 
thickness $2$.
The thickness equals the period of the generating cycle, since both
elements of $\sS^{(0)}$ lie in the interior of the tile $T_1$.
Thus $X_1$ is a substitution Delone set  but not a self-replicating Delone set.
}
\end{exam}

\begin{exam}~\label{ex-selfrep}(Self-replicating Delone set family having
a primitive cycle of order larger than one)
{\rm Let $\he A = [-2]$ and $\sD_{11} = \{-2, -1\}$. Then $\sX = (X_1)$
with $X_1 = \ZZ$ is an irreducible substitution Delone set family
whose fundamental cycle is $\{0, -1\}$ and has period 2. The
corresponding self-affine tile is $T_1 = [0,1]$. So $\sT + \sX = T_1 + X_1$
tiles $\RR$, hence is a self-replicating Delone set.
}
\end{exam}
\begin{exam}~\label{ex75} (Non-periodic and aperiodic
self-replicating Delone set families)
{\rm A Delone set $X$ is (fully) {\em  periodic} if it has
a full rank lattice of periods $\Lambda_X = \{ \bt: X = X+\bt\}$. It
is {\em non-periodic} if it
is not fully periodic and it is
{\em aperiodic} if it has no periods, i.e. $\Lambda_X=\{ {\bf 0} \}.$
An example of a two-dimensional self-replicating
Delone set, which not fully periodic,
but has a one-dimensional lattice of periods, was given in
 Lagarias and Wang~\cite[Example 2.3]{LW96a}. Recently
Lee and Moody~\cite{LM00} constructed many self-replicating
Delone sets which are non-periodic, including aperiodic
examples, whose points are
contained in a lattice in $\RR^d$. They give such examples associated
to non-periodic tilings including the 
sphinx tiling of Godreche~\cite{God89}
and the chair tiling. 
}
\end{exam}


{\tt email: jcl@research.att.com}\\
\hspace*{.85in}{\tt wang@math.gatech.edu}

\end{document}